\documentclass[graybox]{svmult}

\usepackage{type1cm}        % activate if the above 3 fonts are
\usepackage{makeidx}         % allows index generation
\usepackage{graphicx}        % standard LaTeX graphics tool
\usepackage{multicol}        % used for the two-column index
\usepackage[bottom]{footmisc}% places footnotes at page bottom

\usepackage{newtxtext}       % 
\usepackage[varvw]{newtxmath}       % selects Times Roman as basic font
\usepackage{multirow}

\newcommand{\Q}{\mathbf{Q}}
\renewcommand{\u}{\mathbf{u}}

\newcommand{\q}{\mathbf{q}}
\newcommand{\F}{\mathbf{F}}

\renewcommand{\v}{\mathbf{v}}
\newcommand{\B}{\mathbf{B}}
\renewcommand{\S}{\mathbf{S}}
\newcommand{\0}{\mathbf{0}}
\newcommand{\phitilde}{\tilde{\phi}}
\newcommand{\mbf}[1]{\mathbf{#1}}			%
\newcommand{\x}{\mbf{x}}
\newcommand{\BB}[1]{#1}

\makeindex             % used for the subject index
\newcommand{\blue}[1]{{\leavevmode\color{black} #1}}

\begin{document}

\title*{High-order Arbitrary-Lagrangian-Eulerian schemes on \textit{crazy} moving Voronoi meshes}
\author{Elena Gaburro and Simone Chiocchetti}
%\authorrunning{High order ALE schemes  on \textit{crazy} moving Voronoi meshes} 
\institute{Elena Gaburro \at Inria, Univ. Bordeaux, CNRS, Bordeaux INP, IMB, UMR 5251, 200 Avenue de la Vieille Tour, 33405 Talence cedex, France, \email{elena.gaburro@inria.fr}
\and Simone Chiocchetti \at Department of Civil, Environmental and Mechanical Engineering, University of Trento, Via Mesiano, 77 - I-38123 Trento, Italy. \email{simone.chiocchetti@unitn.it}}

\maketitle
% no more than 200 words, repeted twice

\abstract*{
Hyperbolic partial differential equations (PDEs) cover a wide range of interesting phenomena, 
from \textit{human} and \textit{hearth}-sciences up to \textit{astrophysics}:
this unavoidably requires the treatment of many space and time scales
in order to describe at the same time \textit{observer-size} macrostructures, 
\textit{multi-scale} turbulent features, and also \textit{zero-scale} shocks.
Moreover, numerical methods for solving hyperbolic PDEs must reliably handle
different families of waves: 
smooth \textit{rarefactions}, and discontinuities of \textit{shock} and \textit{contact} type. 
In order to achieve these goals, an effective approach consists in the combination of 
space-time-based \textit{high-order} schemes, 
very accurate on smooth features even on coarse grids, 
with \textit{Lagrangian} methods, which, by moving the mesh with the fluid flow, 
yield highly resolved and minimally dissipative results on both shocks and contacts. 
However, ensuring the high quality of \textit{moving meshes} is a 
huge challenge that needs the development of innovative and unconventional techniques.
The scheme proposed here falls into the family of Arbitrary-Lagrangian-Eulerian (ALE) methods,
with the unique additional freedom of evolving the \textit{shape} of the mesh elements  
through connectivity changes. 
We aim here at showing, by simple and very salient examples, the capabilities of high-order ALE schemes, 
and of our novel technique, based on the high-order space-time treatment of \textit{topology changes}.
}

\abstract{
Hyperbolic partial differential equations (PDEs) cover a wide range of interesting phenomena, 
from \textit{human} and \textit{hearth}-sciences up to \textit{astrophysics}:
this unavoidably requires the treatment of many space and time scales
in order to describe at the same time \textit{observer-size} macrostructures, 
\textit{multi-scale} turbulent features, and also \textit{zero-scale} shocks.
Moreover, numerical methods for solving hyperbolic PDEs must reliably handle
different families of waves: 
smooth \textit{rarefactions}, and discontinuities of \textit{shock} and \textit{contact} type. 
In order to achieve these goals, an effective approach consists in the combination of 
space-time-based \textit{high-order} schemes, 
very accurate on smooth features even on coarse grids, 
with \textit{Lagrangian} methods, which, by moving the mesh with the fluid flow, 
yield highly resolved and minimally dissipative results on both shocks and contacts. 
However, ensuring the high quality of \textit{moving meshes} is a 
huge challenge that needs the development of innovative and unconventional techniques.
The scheme proposed here falls into the family of Arbitrary-Lagrangian-Eulerian (ALE) methods,
with the unique additional freedom of evolving the \textit{shape} of the mesh elements  
through connectivity changes. 
We aim here at showing, by simple and very salient examples, the capabilities of high-order ALE schemes, 
and of our novel technique, based on the high-order space-time treatment of \textit{topology changes}.
}

\section{Introduction}
\label{sec.Intro}

In order to reduce as much as possible the numerical errors due to nonlinear convective terms, 
it is possible to exploit the power of Lagrangian methods: 
with this kind of algorithms, the new position and configuration of each element of the mesh is 
recomputed at each timestep according to the local fluid velocity, so that we can closely 
follow the fluid flow in a Lagrangian fashion.
In this way the nonlinear convective terms disappear and Lagrangian schemes 
exhibit \textit{negligible numerical dissipation} at \textit{contact} waves and 
material \textit{interfaces}; 
moreover they results to be \textit{Galilean and rotational invariant}, 
and they provide, without any additional effort, an automatic \textit{mesh refinement} feature 
even when the cell count is maintained constant, simply by transporting the mesh elements wherever needed.

The use of Lagrangian methods dates back to the works of 
\cite{Neumann1950, Wilkins1964} and then many further improvements have been introduced in literature; 
we cite here only some few relevant historical examples and 
review papers~\cite{munz94,CaramanaShashkov1998,LoubereShashkov2004,Despres2009,Maire2010,chengshu2,scovazzi2,Dobrev3,despres2017numerical,morgan2021origins}.

However, ensuring the high quality of a \textit{moving mesh} 
over long simulation times is difficult, therefore a certain degree of flexibility should be
allowed in order to avoid mesh distortion, for example a slightly relaxed choice of the actual mesh
velocity w.r.t the real fluid velocity, as well as the freedom of not only \textit{moving} the
control volumes, but really \textit{evolving} their shapes and allowing topology and neighborhood changes.
This led to the introduction of Arbitrary-Lagrangian-Eulerian~(ALE) schemes 
of \textit{direct}~\cite{ALELTS1D,ALELTS2D,boscheri2013arbitrary} 
and \textit{indirect}~\cite{ReALE2010,ReALE2011,ReALE2015} type.
In particular, as stated and shown in~\cite{springel2010pur}, connectivity changes between 
different time level constitute a valid 
alternative to remeshing~\cite{ReALE2011, ReALE2015, ReALE2010} for preserving or restoring mesh quality
in a Lagrangian setting.

With this in mind, we present here a novel family of very \textit{high-order} direct
Arbitrary-Lagrangian-Eulerian (ALE) Discontinuous Galerkin (DG) schemes for
the solution of general nonlinear hyperbolic PDE systems on moving Voronoi meshes that are
\textit{regenerated} at each timestep and which explicitly allow \textit{topology changes} in time,
in order to benefit simultaneously from high-order methods, high quality grids and substantially
reduced numerical dissipation; this method has been introduced for the first time by the two authors in~\cite{gaburro2020high}.

The \textit{key ingredient} of our approach is the integration of a \textit{space-time conservation
formulation} of the governing PDE system over closed, non-overlapping \textit{space-time} control
volumes~\cite{boscheri2013arbitrary} that are constructed from the moving, \textit{regenerated}, 
Voronoi-type polygonal meshes which are centroid-based dual grids of the Delaunay triangulation of a
set of \textit{generator} points: 
this leads to also consider what we refer to in this work as \textit{crazy degenerate} control 
volumes, or \textit{space-time sliver} elements, that only arise when adopting a space-time framework,
and would not exist from a purely spatial point of view! Such degenerate elements provide a 
clear formal way of handling mesh connectivity changes while preserving the high-order of accuracy of the numerical 
method.

\subsection{Goals}

The goal of this book chapter is to briefly review this novel and promising approach based on high 
order direct ALE schemes with topology changes, and to provide numerical evidence regarding the utility of 
Lagrangian methods in general, and in particular of the new technique object of this work, by means of 
simple and illustrative benchmarks.

\subsection{Structure}

The rest of this chapter is organized as follows. 
We first introduce the class of equations of interest in Section~\ref{sec.pde}. 
Next, we briefly summarize the main characteristics of the employed direct Arbitrary-Lagrangian-Eulerian scheme, 
focusing in particular on the space-time approach and its extension to \textit{crazy} sliver elements, 
whose formation, caused by topology 
changes, will be also addressed in Section~\ref{sec.crazy}.
Then, the core of this work consists in providing numerical evidence 
for \textit{(i)}~the key role of topology changes and sliver elements in a high-order moving mesh code, 
and \textit{(ii)}~the clear advantages of Lagrangian schemes on widely adopted benchmark problems.
Finally, we give some conclusive remarks and an outlook towards future work in Section~\ref{sec.conclusion}.

\section{Hyperbolic partial differential equations}
\label{sec.pde}

In order to model a wide class of physical phenomena, we consider 
a very general formulation of the governing equations, namely all those which can be described by 
\begin{equation} 
\label{eq.generalform}
\partial_t{\Q} + \nabla \cdot \F(\Q) + \B(\Q) \cdot \nabla \Q = \mathbf{S}(\Q),
\end{equation}
where $\Q$ is the vector of the conserved variables, $\F$ the non linear flux, 
$\B\cdot \nabla \Q$ the nonconservative products, and $\mathbf{S}$ a nonlinear algebraic source term. 
Many physical models can be cast in this form, from the simple shallow water 
system, some multiphase flow models, the magnetohydrodynamics equations, 
up to the Einstein field equations of general relativity (with appropriate reformulation) or the GPR 
unified model of continuum mechanics, see for 
example~\cite{fambri2017space,gaburro2018wellb,kemm2020simple,gaburro2021well,chiocchetti2021high,tavelli2020space,chiocchetti2020solver,chiocchetti2021high,gabriel2021unified,peshkov2021simulation,boscheri2022cell}; 
in this work, we will present illustrative results concerning the Euler equations of gasdynamics.

\section{Numerical method}
\label{sec.method}

In this Section we presents a concise description of our direct Arbitrary-Lagrangian-Eulerian (ALE) 
Discontinuous Galerkin (DG) scheme on moving Voronoi-type meshes with topology changes; for any 
additional details we refer to our recent paper~\cite{gaburro2020high}.

\medskip 

At the beginning of the simulation, we discretize our moving domain by a centroid-based Voronoi-type
tessellation built from a set of generators (the orange points in Figure~\ref{fig.controlvolumes}), 
and we represent our data, the conserved variables $\Q$,
via discontinuous high-order polynomials in each mesh polygon (we 
indicate the degree of the polynomial representation by $P_N$).
Then, we let the generators move
with a velocity chosen to be as close as possible to the local fluid velocity, computed mainly from
a high-order approximation of their pure Lagrangian trajectories, with small corrections obtained
from a flow-adaptive mesh optimization technique. The positions of the generators are being
continuously updated, and thus their Delaunay triangulation may change at any timestep
and the same will hold for the dual polygonal tessellation.
Then, a space-time connection between two polygonal tessellations corresponding to two successive time levels
has to be established in order to
evolve the solution in time locally and integrate the governing PDE. 

\subsection{Direct Arbitrary-Lagrangian-Eulerian schemes} 

\begin{figure}[!bp]
	{\includegraphics[width=0.33\linewidth]{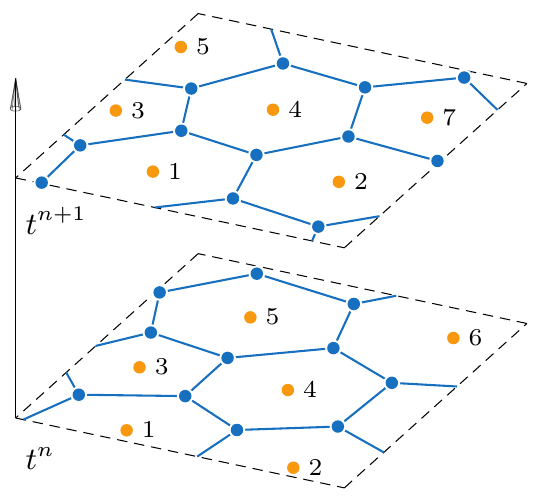}}%
	{\includegraphics[width=0.33\linewidth]{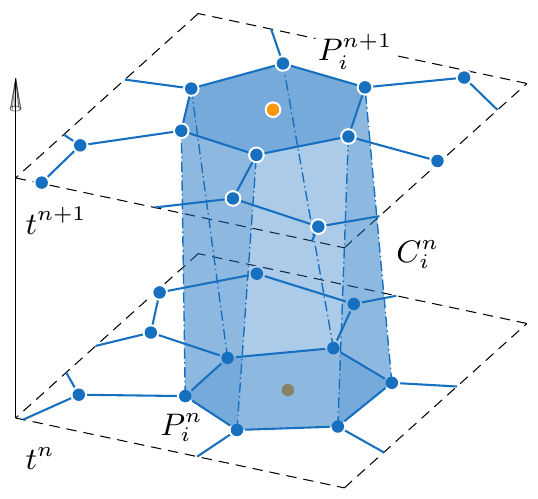}}%
	{\includegraphics[width=0.33\linewidth]{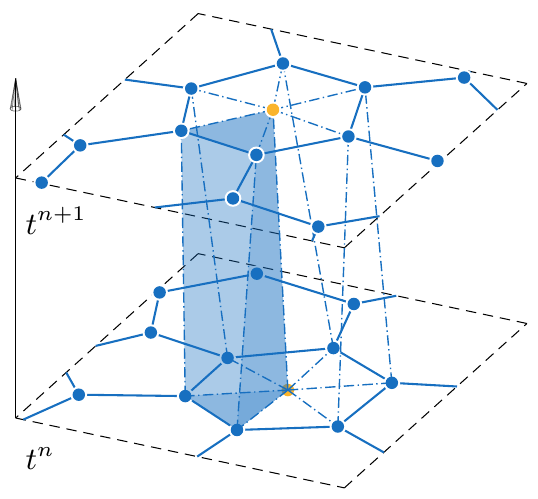}}%
	\caption{Space time connectivity \textit{without} topology changes, main space-time control volume (middle) and a standard sub-space-time control volume (right).}
	\label{fig.controlvolumes}
\end{figure}

The key idea of \textit{direct} ALE methods (in contrast to
\textit{indirect} ones) consists in connecting two tessellations by means of so-called
\textit{space-time control volumes} $C_i^n$, and recover the unknown solution at the new 
timestep $\u_h^{n+1}$ \textit{directly} inside the new polygon $P_i^{n+1}$,
from the data available at the previous timestep $\u_h^{n}$ in $P_i^n$. 
This is achieved through the \textit{integration}, over such control volumes, of the fluxes, the nonconservative 
products and the source terms, by means of a high-order fully discrete predictor-corrector ADER 
method~\cite{dumbser2008unified,gaburro2021posteriori}. 
In this way, the need for any further remapping/remeshing steps is 
totally eliminated. By adopting the tilde symbol for referring to space-time quantities, 
our direct ALE scheme~\cite{gaburro2020high,gaburro2021unified} reads
\begin{equation}
\label{eqn.ALE-ADER}
\begin{aligned}
	\int_{P_i^{n+1}} 
	\phitilde_k  \u_h^{n+1}  &=
	\int_{P_{i}^n} 
	\phitilde_k  \u_h^n   
	- \sum_{j} \int_{\partial C_{ij}^n} 
	\phitilde_k  \mathcal{F}(\q_h^{n,-},\q_h^{n,+}) \cdot \mathbf{\tilde n}  
	+ \int_{C_i^n} 
	\tilde \nabla \phitilde_k \cdot \tilde{\F} (\q_h^n) \\
	&+ \int_{C_i^n} \phitilde_k \left ( \mathbf{S}(\q_h^n) - \tilde{\B}(\q_h^n)\cdot \nabla \q_h^n \right), 
\end{aligned}
\end{equation}
where $\phitilde_k$ is a set of moving space-time basis functions, while $\q_h^{n,+}$ and $\q_h^{n,-}$ are 
high-order space-time extrapolated data computed through the ADER predictor. 
Finally, $\mathcal{F}(\q_h^{n,-},\q_h^{n,+})$ is an ALE numerical flux function which takes into 
account fluxes across space-time cell boundaries $\partial C_{ij}^n$ as well as jump terms 
related to nonconservative products.
In particular, \blue{we adopt here a  two-point path-conservative numerical flux function of Rusanov-type~\cite{Rusanov:1961a,pares2006numerical}}
\begin{equation}
\begin{aligned} 
	\mathcal{F}(\q_h^{n,-},\q_h^{n,+}) \cdot \mathbf{\tilde n} = &  
	\frac{1}{2} \left( \tilde{\F}(\q_h^{n,+}) + \tilde{\F}(\q_h^{n,-})  \right) \cdot \mathbf{\tilde n}_{ij}  - 
	\frac{1}{2} s_{\max} \left( \q_h^{n,+} - \q_h^{n,-} \right) \\
	&+ \frac{1}{2} \left(\int \limits_{0}^{1} \tilde{\mbf{B}} \left(\mbf{\Psi}(\q_h^{n,-},\q_h^{n,+},s) \right)\cdot\mbf{n} \, d\x \right)
	\cdot\left(\q_h^{n,+} - \q_h^{n,-}\right), 
	\label{eq.rusanov} 
\end{aligned} 
\end{equation} 
where $s_{\max}$ is the maximum eigenvalue of the \textit{ALE Jacobian matrices} evaluated 
on the left and right of the space-time interface and the path $\mbf{\Psi}=\mbf{\Psi}(\q_h^-,\q_h^+, s)$ is a straight-line segment path
connecting $\q_h^{n,-}$ and $\q_h^{n,+}$.

We emphasize that the ALE Jacobian matrix is obtained by subtracting the local normal mesh velocity 
from the diagonal entries 
of the system matrix of the quasilinear form of the governing equations~\cite{toro2013riemann} 
(the Jacobian of the interface-normal flux for conservative systems), thus, when the mesh 
velocity is sufficiently close to the local fluid velocity, 
the wavespeed estimates obtained from the eigenvalues are significantly reduced, 
leading to a lower associated numerical dissipation than what would be mandated in the Eulerian context.
This, especially but not exclusively, in conjunction with complete approximate Riemann solvers~\cite{hllem},
explains the capability of tracking material interfaces and capturing contact 
discontinuities \blue{which are characteristic} of Lagrangian-type schemes.

\begin{figure}[!bp]
	\includegraphics[width=0.33\linewidth]{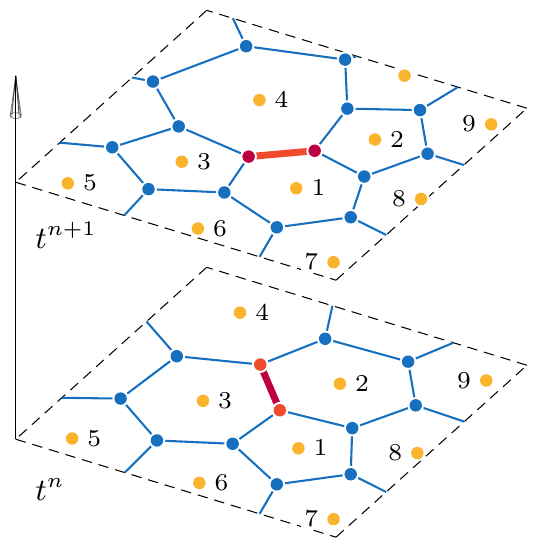}%
	\includegraphics[width=0.33\linewidth]{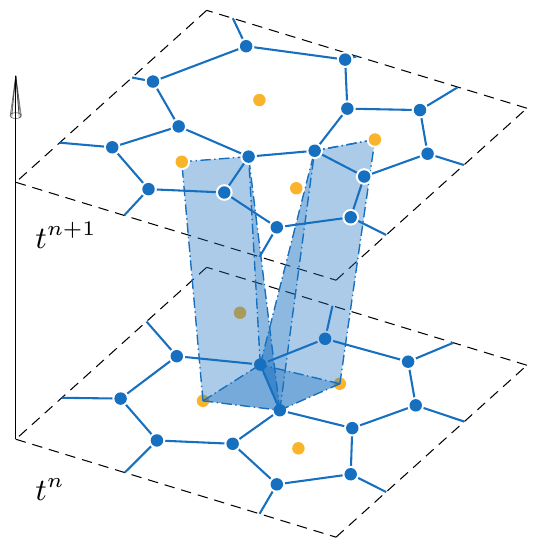}%
	\includegraphics[width=0.33\linewidth]{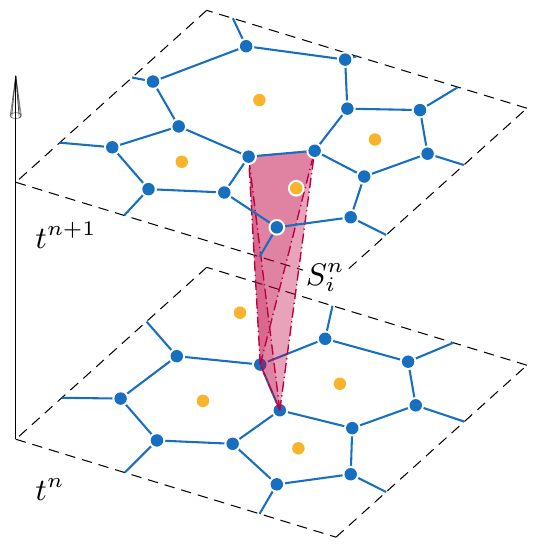}%
	\caption{Space time connectivity \textit{with} topology changes, degenerate sub-space-time control volumes (middle) and \textit{crazy} sliver element (right).}
	\label{fig.crazy}
\end{figure}

Next, in order to compute the integrals with high order of accuracy, complete knowledge of the
\textit{space-time connectivity} between two consecutive timesteps is required, 
as opposed to only the \textit{spatial} information at the two time levels, 
which would be enough for a low order scheme~\cite{springel2010pur}
or for indirect schemes~\cite{ReALE2010,ReALE2011,ReALE2015}.

\begin{figure}[!bp]
	{\includegraphics[width=0.33\linewidth]{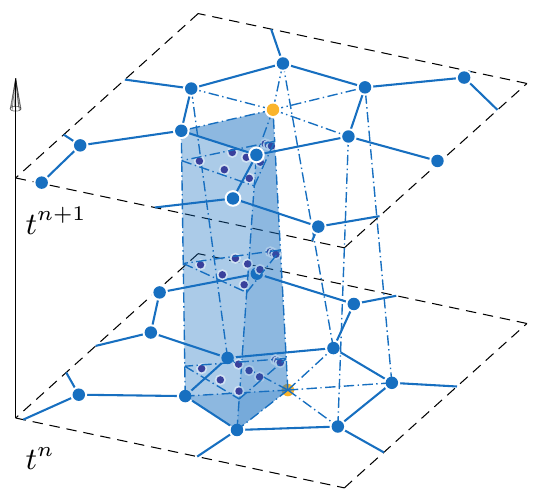}}%
	{\includegraphics[width=0.33\linewidth]{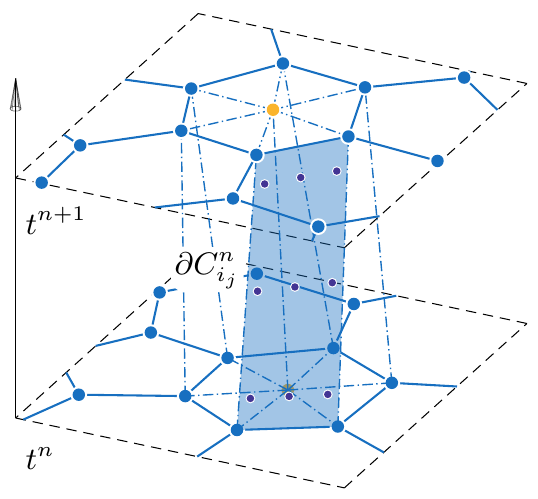}}%
	{\includegraphics[width=0.33\linewidth]{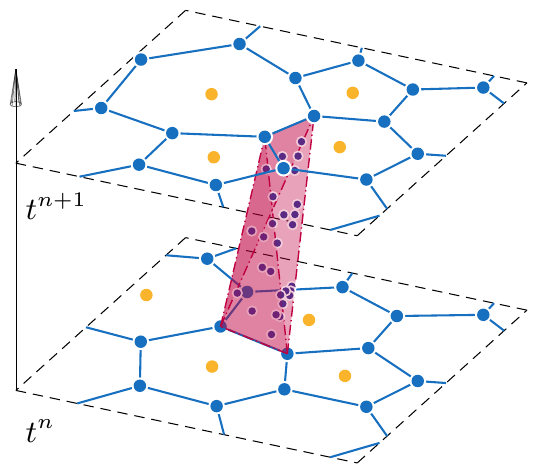}}%
	\caption{Space-time quadrature points for third order methods on standard elements (left), lateral faces (middle) and \textit{crazy} sliver elements (right).}
	\label{fig.quadraturepoints}
\end{figure}

When no topology changes occur, the space-time
geometrical information is easily constructed by connecting via straight line segments the
corresponding vertexes of each polygon, obtaining an oblique prism than can be further subdivided 
into a set of triangular oblique sub-prisms on which quadrature points are readily available (see
Figures~\ref{fig.controlvolumes} and~\ref{fig.quadraturepoints}).

\subsection{Topology changes and \textit{crazy} sliver elements}
\label{sec.crazy}

On the contrary, when a topology change occurs, as in Figure~\ref{fig.crazy}, i.e. the number of
edges, the shape, and the neighbors of a polygon evolve within two consecutive timesteps, the
space-time connection between the mesh elements gives raise to degenerate elements of two types:
\textit{(i)}~\textit{degenerate} \textit{sub}-space-time control volumes, where either the top or bottom
faces are degenerate triangles that are collapsed to a segment; \textit{(ii)}~and also
\textit{crazy sliver} space-time elements $S_i^n$. 
The first type of degenerate elements does not
pose any problems, and was already treated in~\cite{gaburro2017direct}. 
Instead, space-time sliver elements are a completely new type of control volume. 
In particular, they do not exist neither at time $t^n$, nor
at time $t^{n+1}$, since they coincide with an edge of the tessellation at the old \textit{and} at the new
time levels, and for this reason have zero area in space at the two bounding time levels. 
However, they have a \textit{non-negligible volume} in space-time. 
The difficulties related to this kind of elements are due to the fact that
for them an initial condition is not clearly defined at time $t^n$, and that contributions across
these elements should not be lost at time $t^{n+1}$, in order to ensure conservation. All the
details on how to successfully extend our direct ALE scheme also to \textit{crazy} elements can be found in
our recent paper~\cite{gaburro2020high}.

We would like to emphasize that topology changes are fundamental for long time simulations in the
ALE framework, in order to avoid explicit data remap steps, and our \textit{crazy} sliver
elements represent a \textit{novel and formally grounded}
way to allow for a relatively simple space-time connection around a change of connectivity.
The numerical results shown in Section~\ref{sec.shu} provide a clear proof of the necessity of 
topology changes already on the simple situation of the long time evolution of a stationary isentropic smooth vortex.

\subsection{ADER-ALE algorithm: the predictor step} 

The \textit{predictor} step represents an essential ingredient
for obtaining high-order in time in a fully-discrete one-step procedure: it yields a \textit{local} 
solution of the governing equations~\eqref{eq.generalform} \textit{in the small} $\q_h^n$, inside each
space-time element, including the \textit{crazy} elements. The solution is \textit{local} in the sense that
it is obtained by only considering the initial data in each polygon, the governing equations and the
geometry of $C_i^n$, without taking into account interactions between $C_i^n$ and its neighbors.
Such local solution is computed for each standard space-time control volume $C_i^n$ 
and for each \textit{crazy} control volume $S_i^n$, in the form of a high
order polynomial \textit{in space and in time}, which serves as a predictor solution, 
to be used for evaluating all the integrals in the \textit{corrector} step~\eqref{eqn.ALE-ADER}, i.e. the
final update of the solution from $t^n$ to $t^{n+1}$. 

\subsection{\textit{A posteriori} sub-cell FV limiter} 

High-order schemes that can be seen as linear in the sense of Godunov~\cite{godunov}, may develop
spurious oscillations in presence of discontinuities. In order to prevent this phenomenon, in the
case of a DG discretization we adopt an \textit{a~posteriori} limiting procedure based on the MOOD paradigm
\cite{CDL1,ADERMOOD,gaburro2021posteriori}: we first apply our unlimited ALE-DG scheme everywhere, and then (\textit{a
posteriori}), at the end of each timestep, we check the reliability of the obtained solution in
each cell against physical and numerical admissibility criteria, such
as floating point exceptions, violation of positivity or violation of a relaxed discrete maximum 
principle (and see~\cite{guermond2018second,KENAMOND2021110254} for further criteria).
Next, we mark as \textit{troubled} those cells where the DG solution cannot be accepted. 
For the troubled cells we now repeat the time evolution by employing, instead of the DG scheme, 
a more robust finite volume (FV) method. Moreover, in
order to maintain the accurate resolution of our original high-order DG scheme, which would be lost
when switching to a FV scheme, the FV scheme is applied on a \textit{finer sub-cell grid}
that accounts for recovering the optimal accuracy of the numerical method performing a reconstruction step.

\section{Numerical examples}
\label{sec.examples}

In order to provide simple and clear numerical evidence of the effectiveness of the proposed ALE scheme with 
topology changes we consider here the well known Euler equations, that can be cast in the form~\eqref{eq.generalform} by choosing 
\begin{equation}
	\label{eulerTerms}
	\Q = \left( \begin{array}{c} \rho   \\ \rho u  \\ \rho v  \\ \rho E \end{array} \right), \quad
	\mathbf{F} = \left( \begin{array}{ccc}  \rho u       & \rho v        \\ 
		\rho u^2 + p & \rho u v          \\
		\rho u v     & \rho v^2 + p      \\ 
		u(\rho E + p) & v(\rho E + p)   
	\end{array} \right), \qquad \B = \0, \qquad \S = \0.  
\end{equation}
The vector of conserved variables $\Q$ is composed of the fluid density $\rho$, 
the momentum density vector $\rho \v=(\rho u, \rho v)$ and the total energy density $\rho E$; 
next, the fluid pressure $p$ is computed using the equation of state for an ideal gas    
\begin{equation}
	\label{eqn.eos} 
	p = (\gamma-1)\left(\rho E - \frac{1}{2} \rho \mathbf{v}^2 \right), 
\end{equation}
where $\gamma$ (in this work taken to be $\gamma = 7/5$) is the ratio of specific heats. For this choice
of equation of state, the adiabatic
speed of sound takes the form $c=\sqrt{{\gamma p}/{\rho}}$. 

\medskip

In what follows we will present numerical results regarding the following notable features of Lagrangian schemes
and of our direct Arbitrary-Lagrangian-Eulerian method with variable topology:
\begin{description}
	\item[{i.}] Flows characterized by strong differential rotations, for example vortices, 
	can be studied over very long periods only by conceding to the element motion the 
	additional freedom of introducing topology changes, see Section~\ref{sec.shu};
	\item[{ii.}] The use of sliver elements allows to clearly define the space-time evolution of 
	the solutions in-between discrete time levels and achieves high-order of accuracy also 
	in presence of many topology changes, see Section~\ref{sec.order};
	\item[{iii.}] Lagrangian schemes sharply capture \textit{shock} 
	waves thanks to the automatic refinement obtained at the shock locations without 
	needing to increment the number of mesh elements but simply because 
	the element density increases wherever needed, see Section~\ref{sec.sedov};
	\item[{iv.}] Lagrangian schemes minimize dissipation of \textit{contact} discontinuities, by applying
	reduced numerical dissipation when using approximate Riemann solvers. In a pure Lagrangian
	context, schemes capable of capturing stationary discontinuities exactly will do the same 
	also for moving interfaces (since the mesh motion is specified to follow such features). 
	Moreover, even when such hard constraints are relaxed in 
	Arbitrary-Lagrangian-Eulerian methods and even using 
	simpler solvers like the Rusanov flux, 
	the bounding wavespeed estimates and the associated numerical dissipation can be much lower 
	than what would be mandated in the Eulerian context,
	% $\lambda = \mathbf{v}\cdot\mathbf{n} \pm c$,
	see Section~\ref{sec.sod};
	\item[{v.}] Lagrangian schemes discretely preserve the Galilean and rotational invariance properties
	of the governing equations, 
	so that they can better capture any events (such as the explosion-type problems reported in this work) 
	that may occur in superposition to a high-speed background flow, see Section~\ref{sec.sod}.
\end{description}

\subsection{Long time evolution of a Shu-type vortical equilibrium}
\label{sec.shu}

\begin{figure}[!bp]  % 3 colonne
	\centering
	\includegraphics[width=0.33\linewidth]{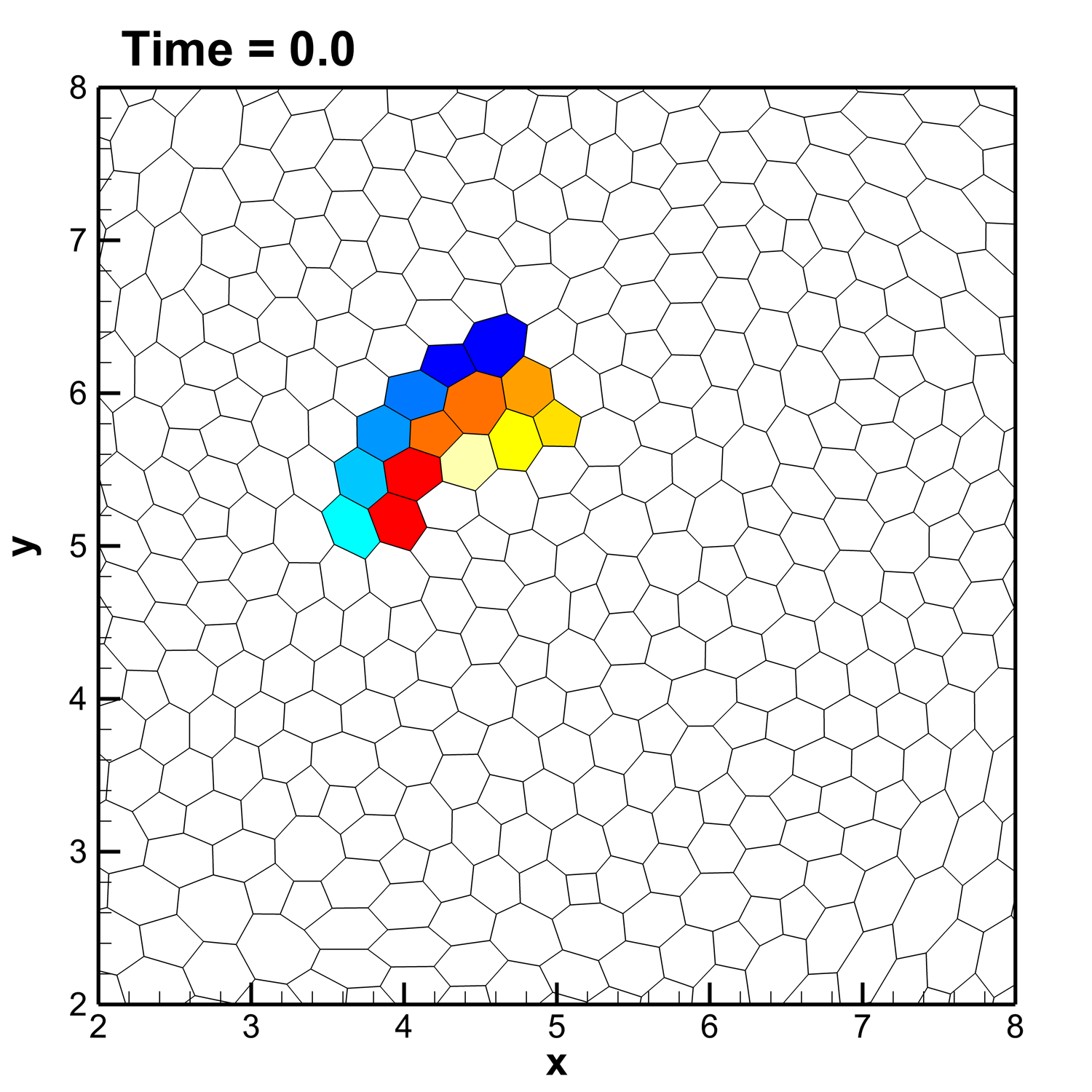}%
	\includegraphics[width=0.33\linewidth]{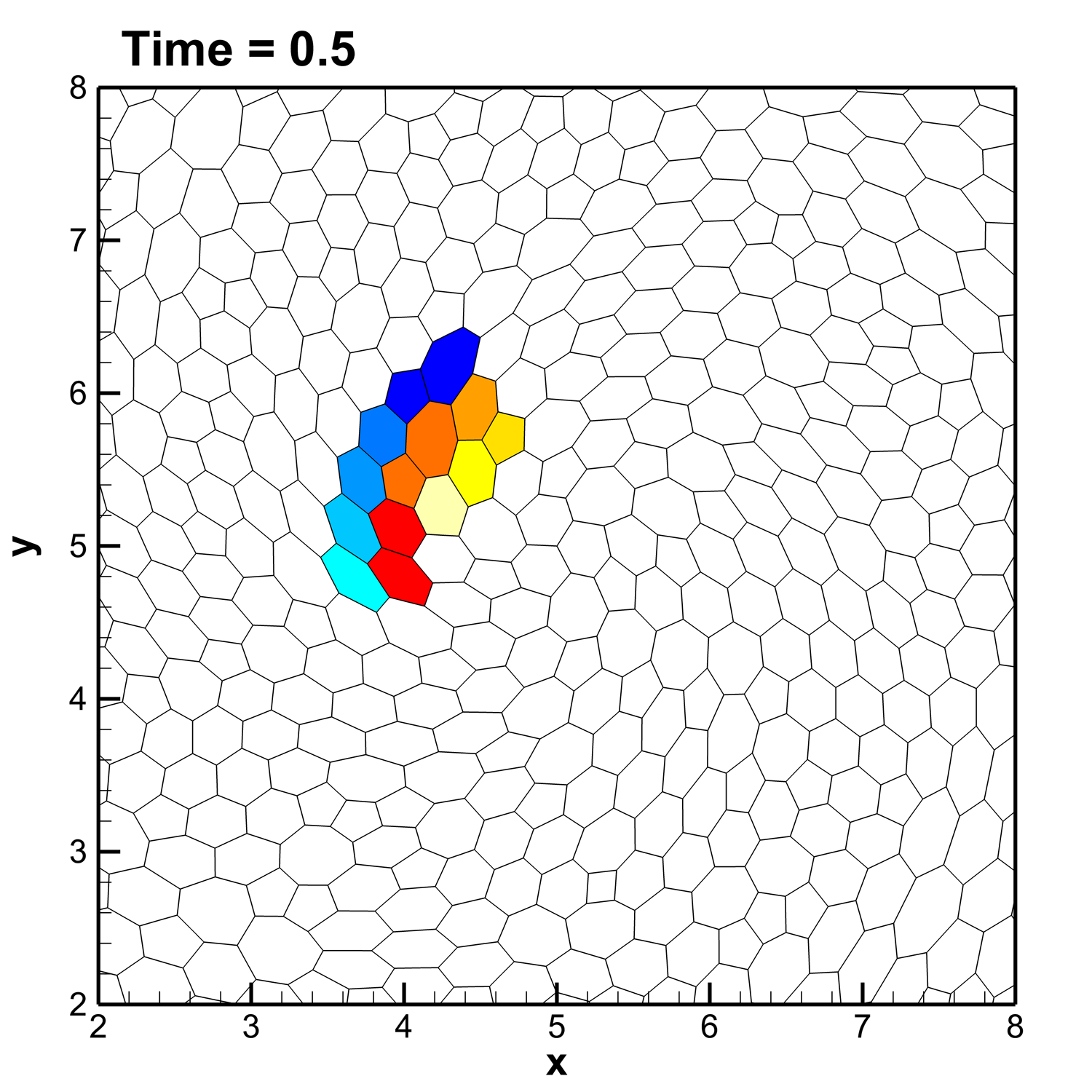}%
	\includegraphics[width=0.33\linewidth]{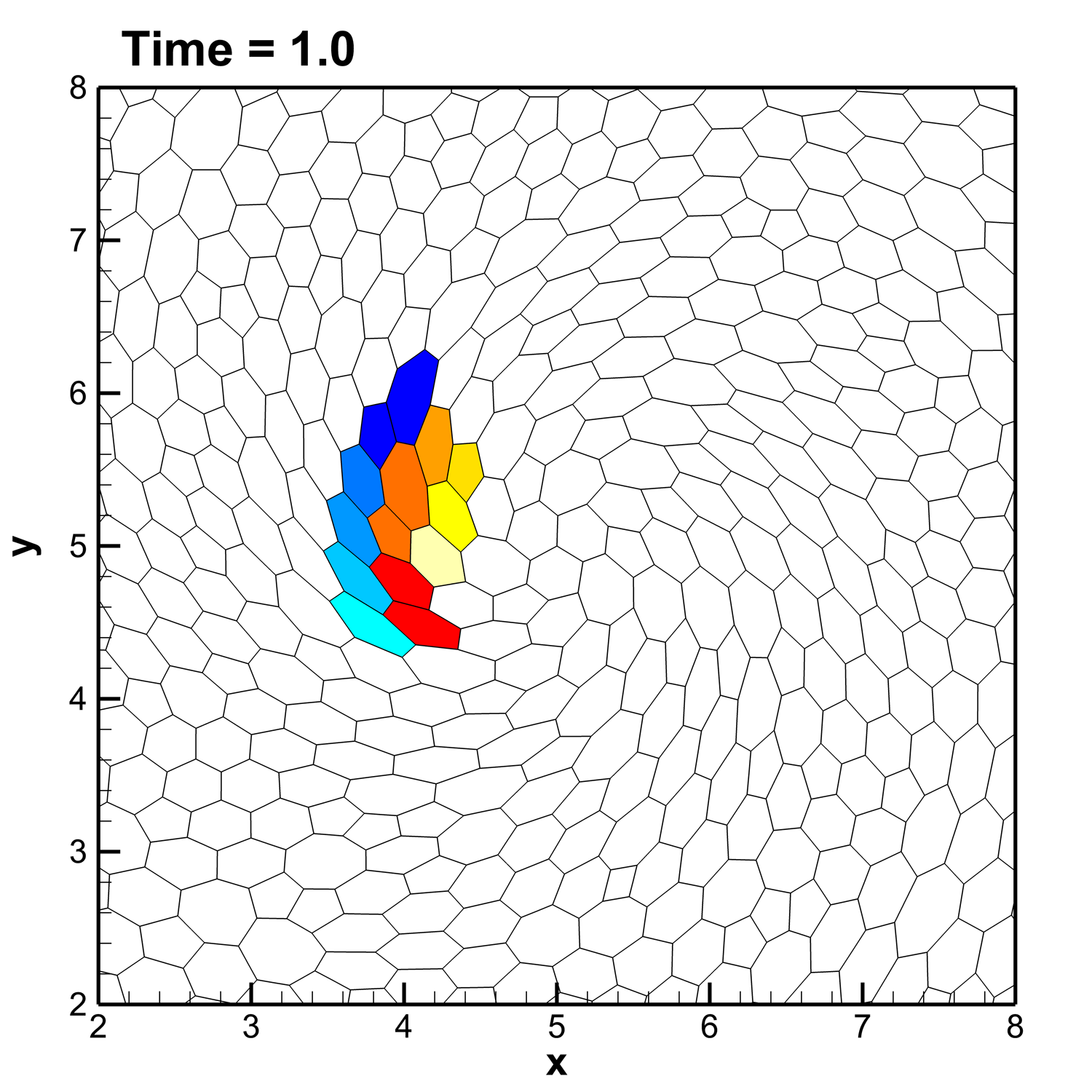}\\
	\includegraphics[width=0.33\linewidth]{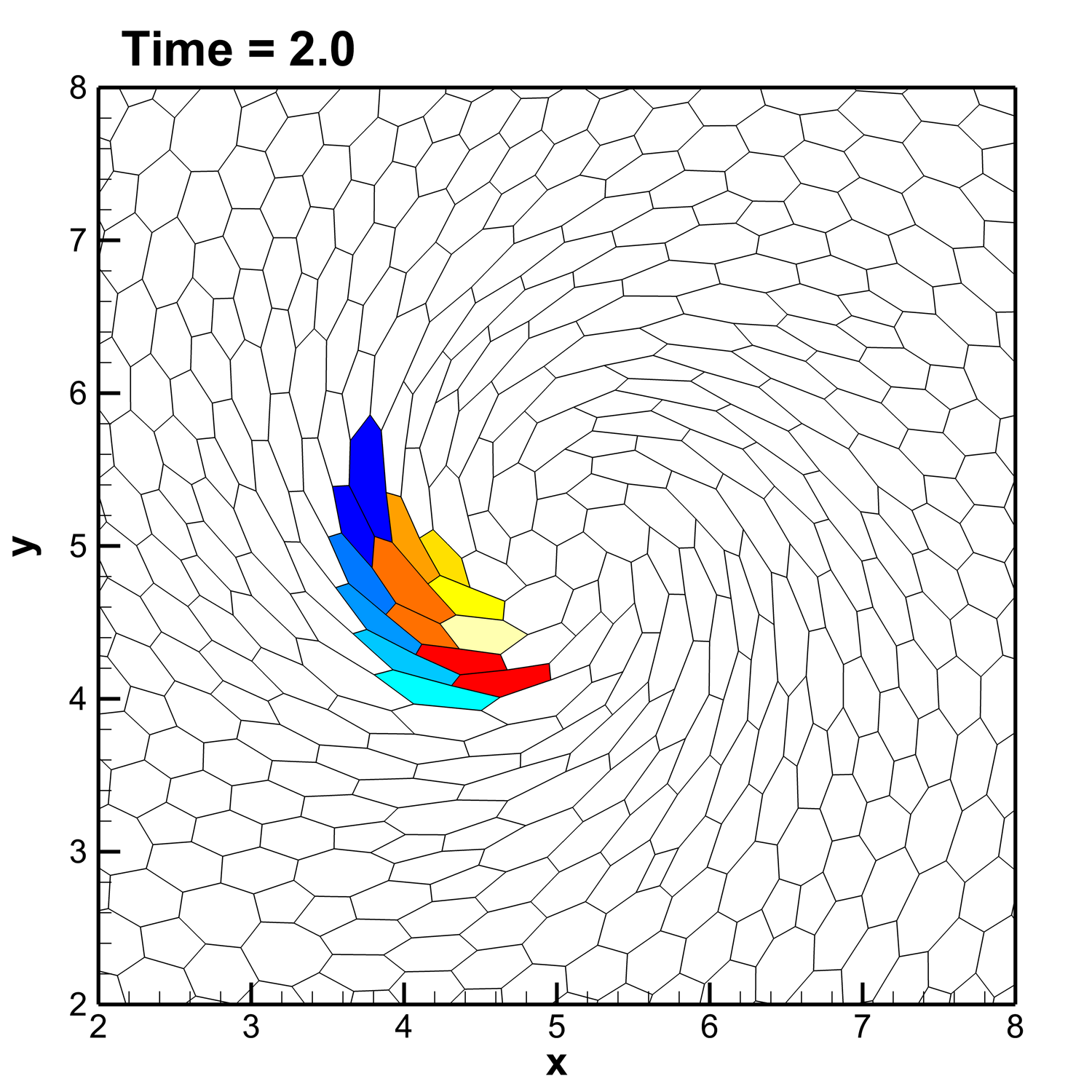}%
	\includegraphics[width=0.33\linewidth]{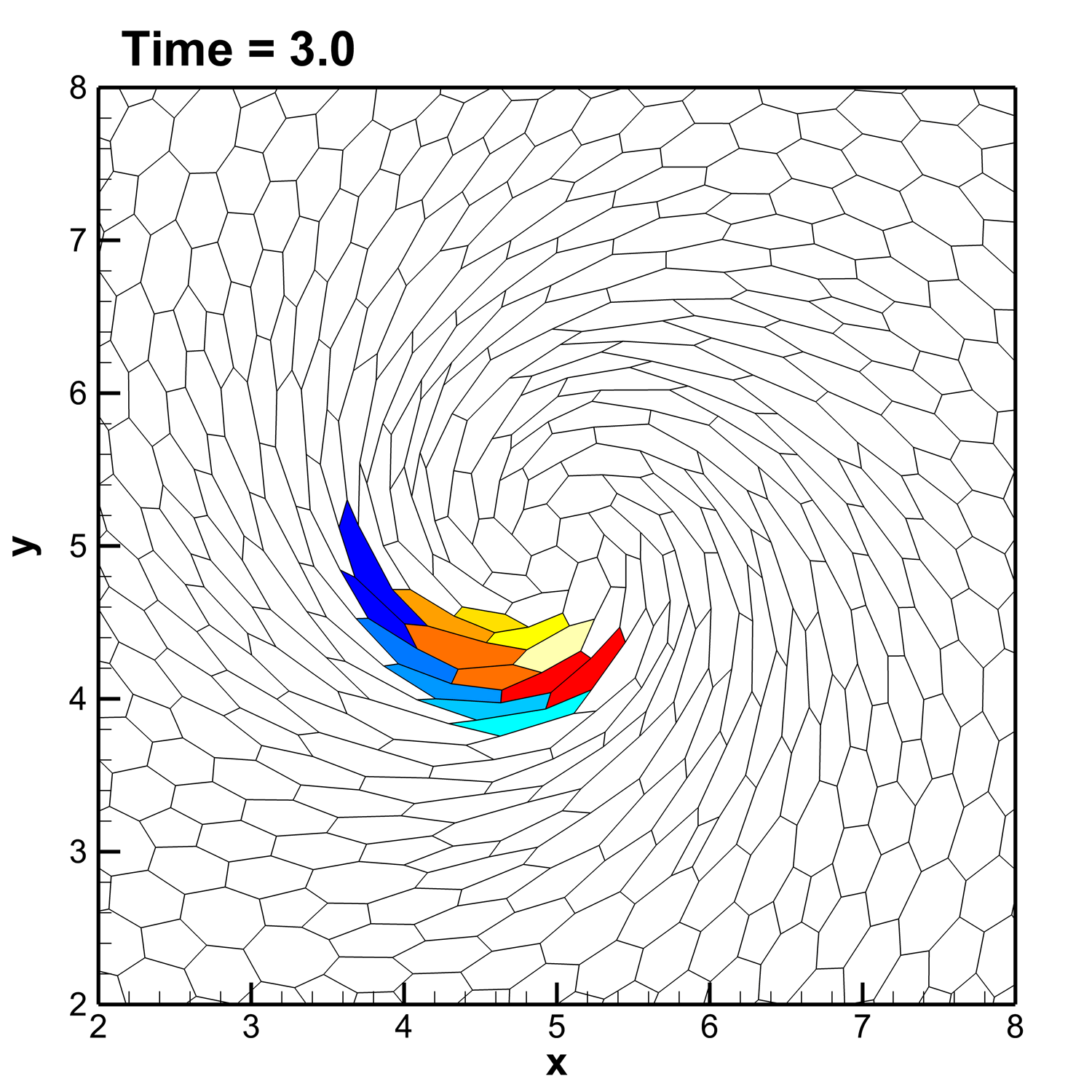}%
	\includegraphics[width=0.33\linewidth]{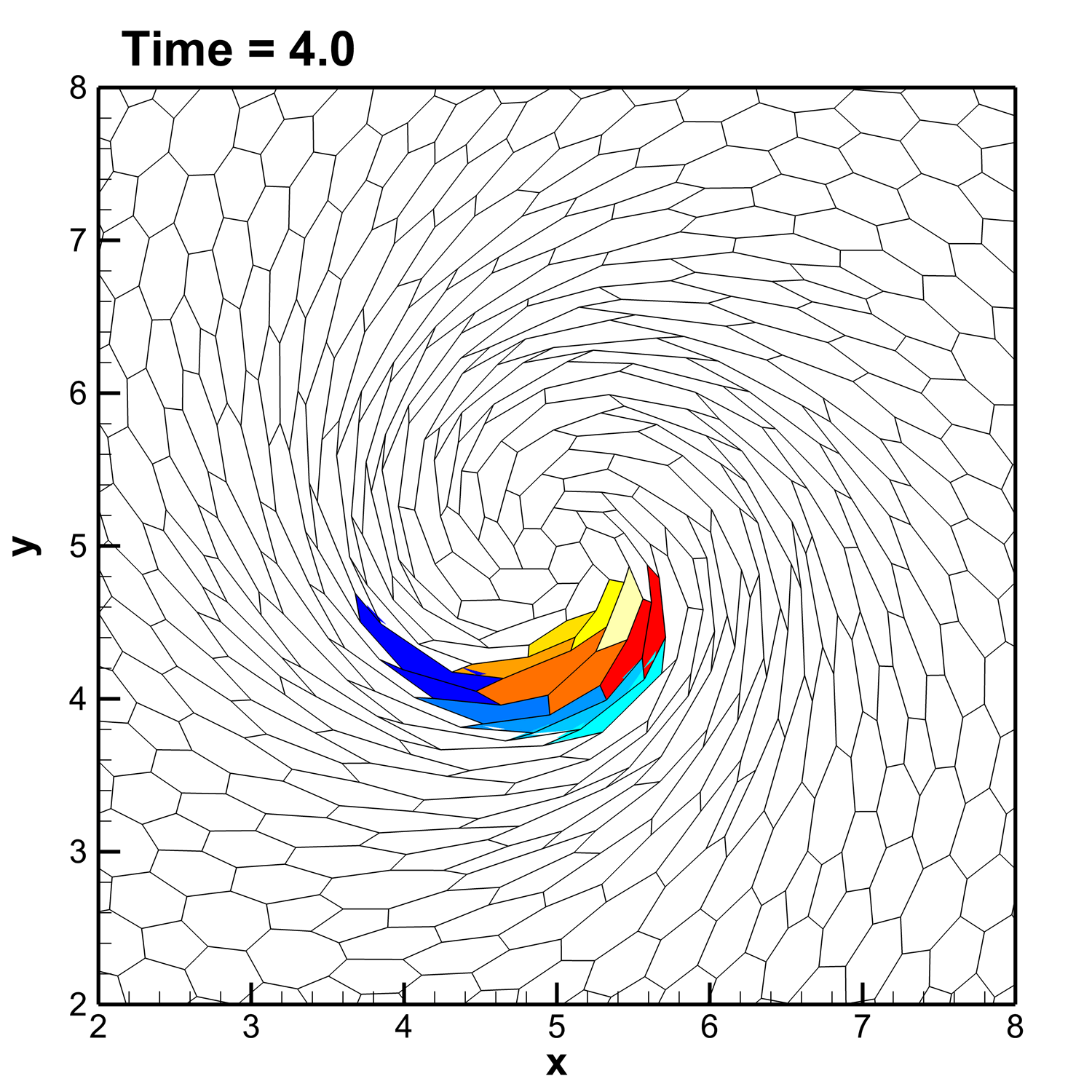}\\
	\caption{Mesh evolution corresponding to the solution of the stationary rotating vortex of Section~\ref{sec.shu} 
		solved on a moving grid with \textit{fixed} topology. The mesh quality rapidly deteriorates: elements 
		are stretched, the timestep size is reduced, and even mesh-tangling occurs, which means that the simulation 
		may stop entirely.}
	\label{fig.shu_rho_ie_notopc}
\end{figure}

\begin{figure}[!bp]  % 3 colonne
	\centering
	\includegraphics[width=0.33\linewidth]{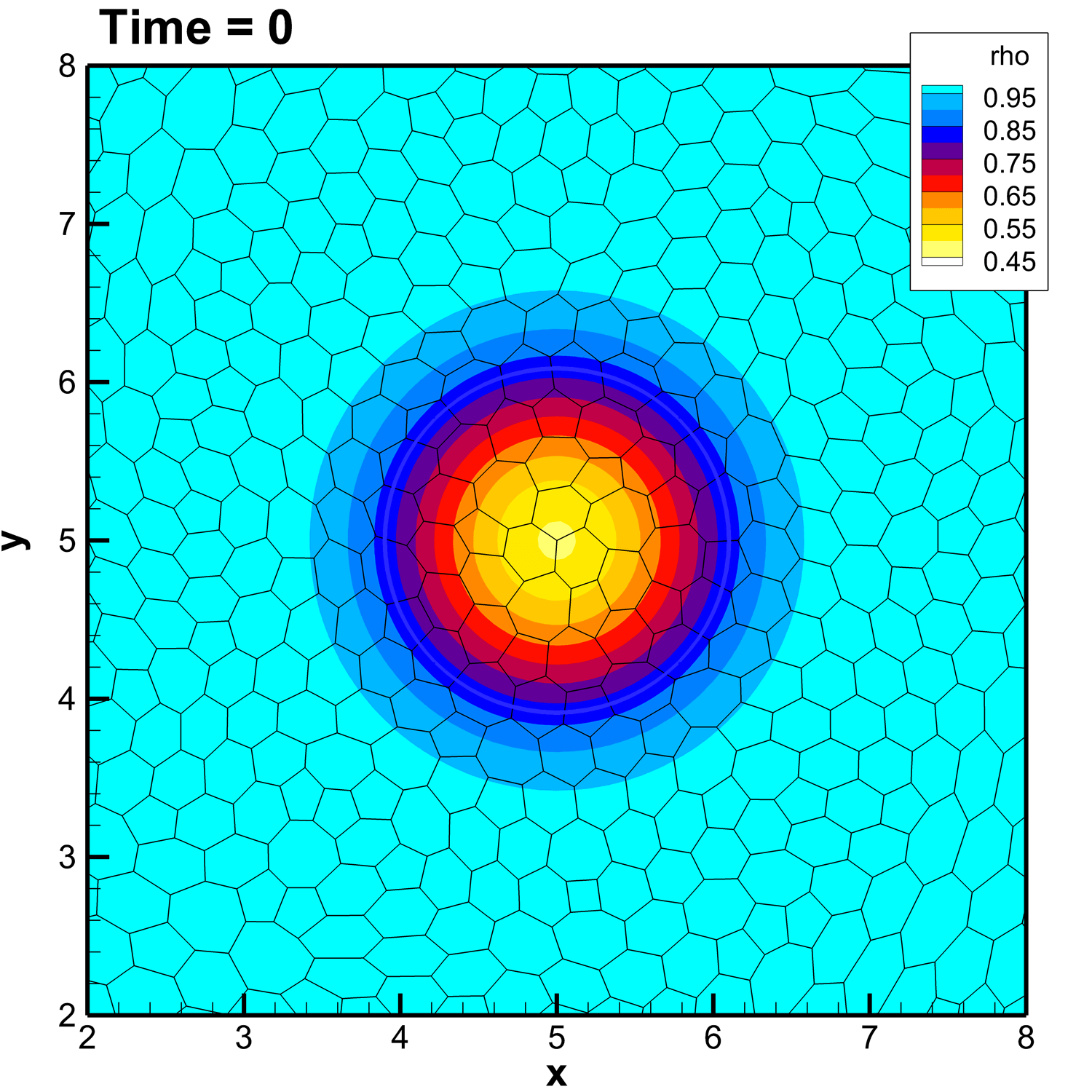}%
	\includegraphics[width=0.33\linewidth]{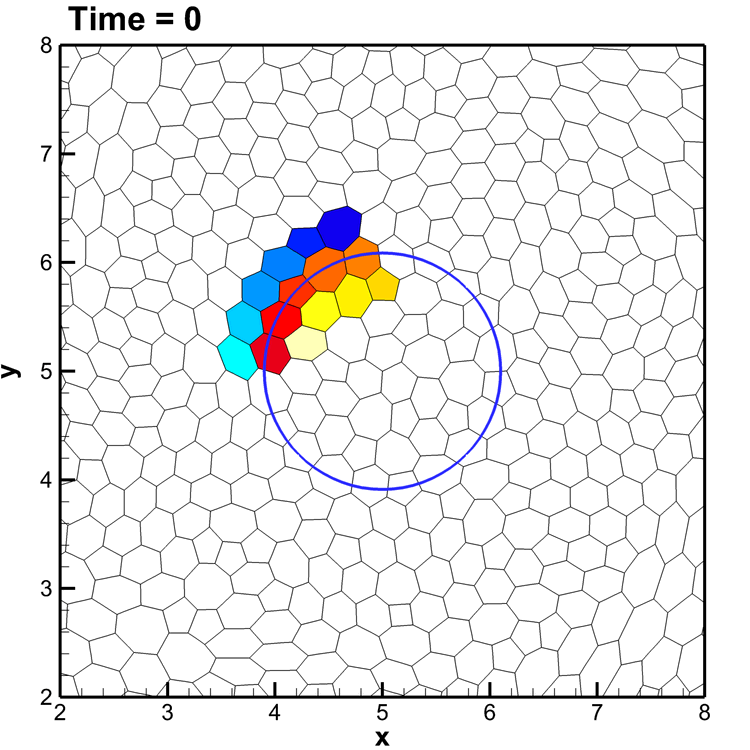}%
	\includegraphics[width=0.33\linewidth]{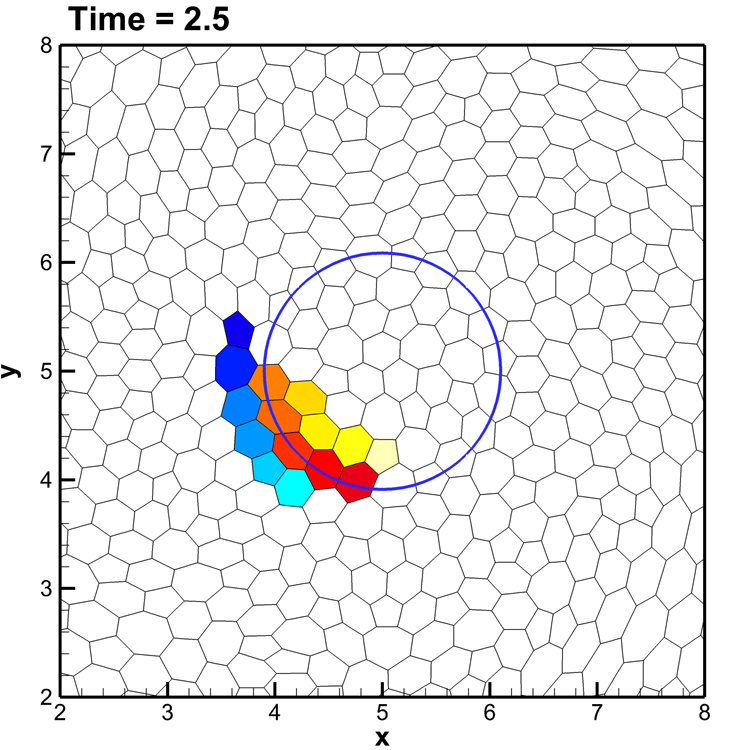}\\
	\includegraphics[width=0.33\linewidth]{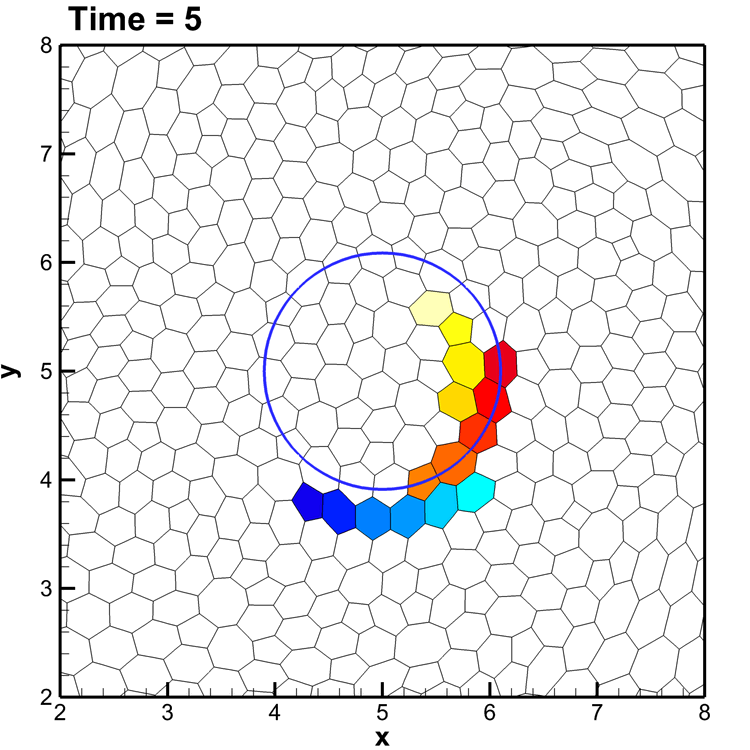}%
	\includegraphics[width=0.33\linewidth]{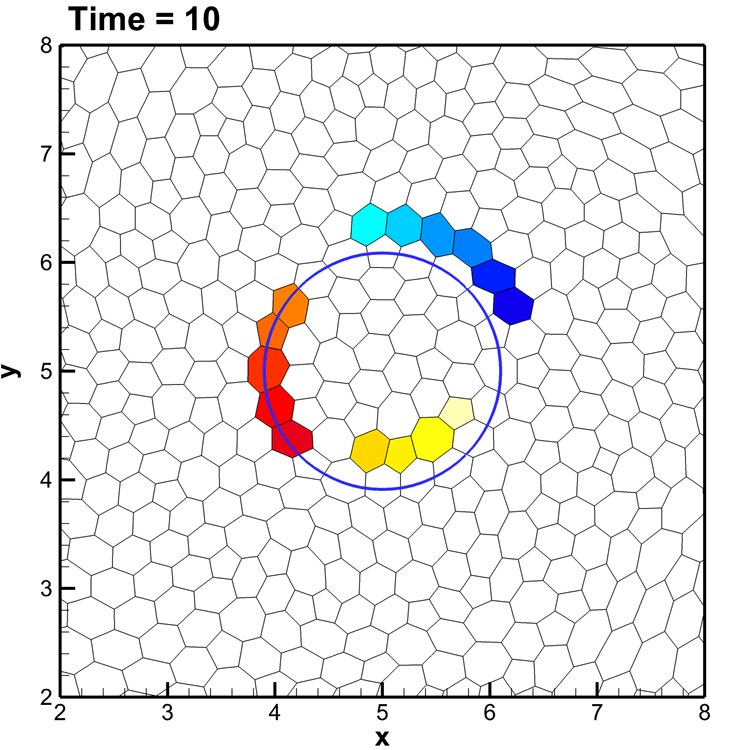}%
	\includegraphics[width=0.33\linewidth]{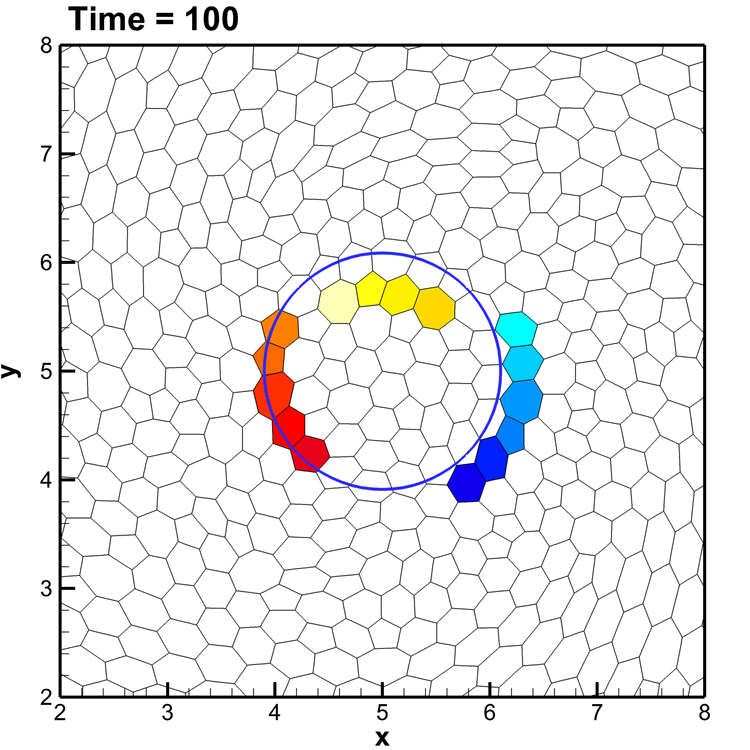}\\
	\includegraphics[width=0.33\linewidth]{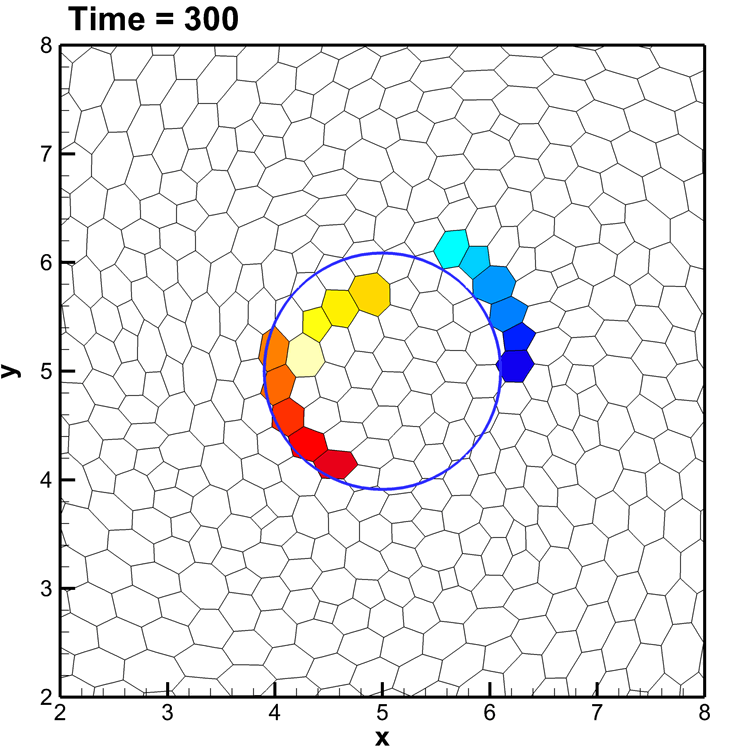}%
	\includegraphics[width=0.33\linewidth]{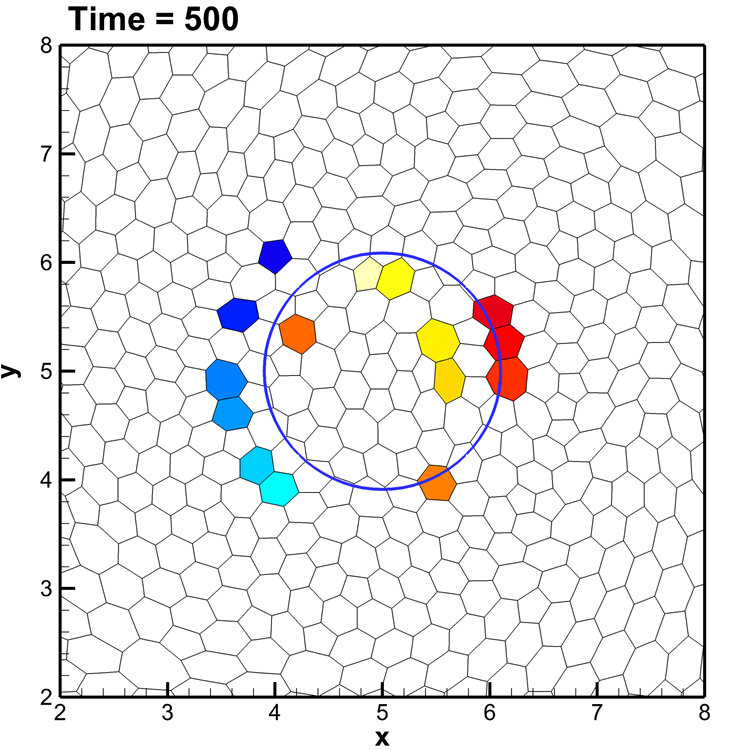}%
	\includegraphics[width=0.33\linewidth]{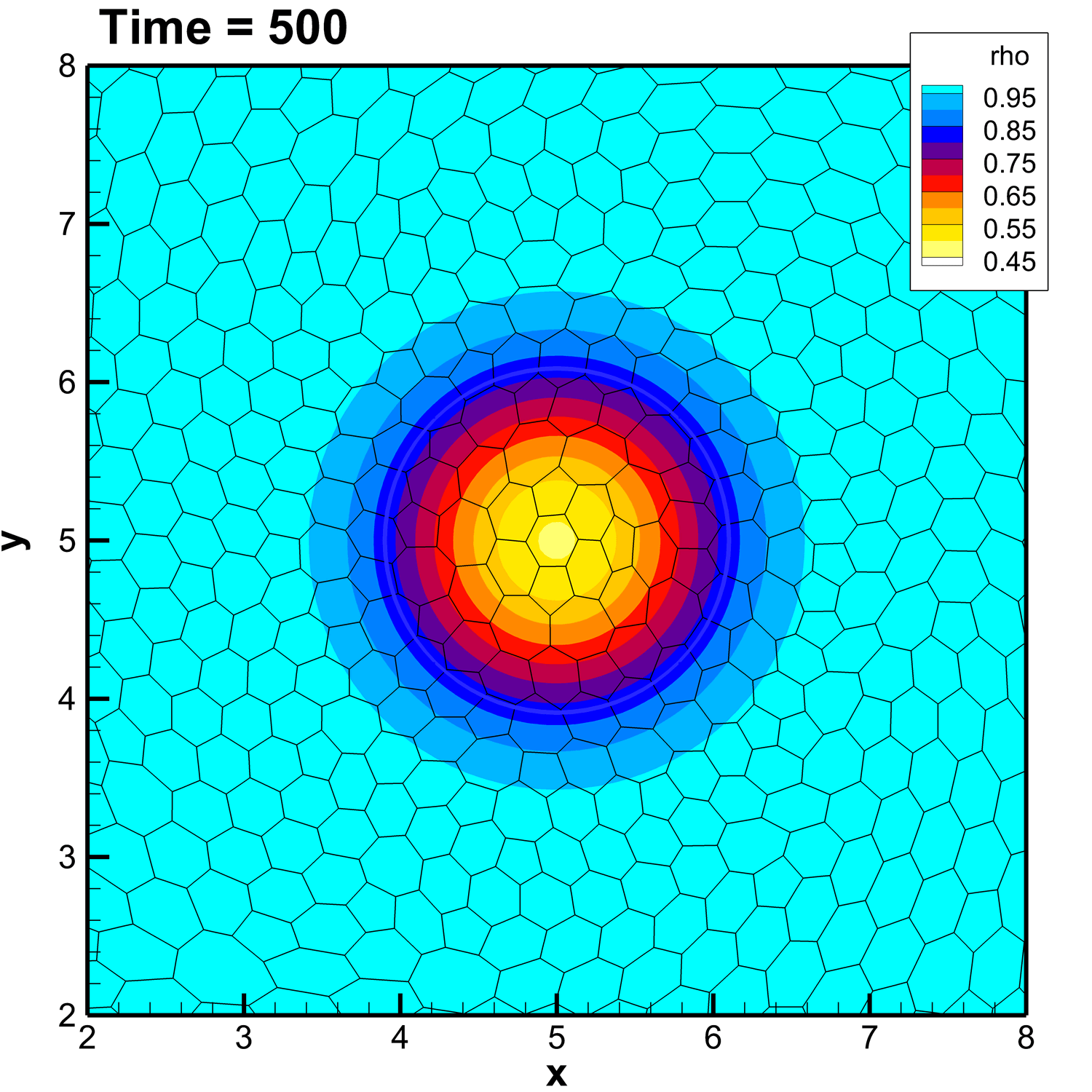}%
	\caption{Stationary rotating vortex solved with our fourth order ALE-DG scheme. 
		Density contours at $t=0$ and $t=500$ and position of a bunch of highlighted 
		elements at different times. Note that the solution is well preserved for more than 
		\textit{eighty} complete rotation periods of the yellow elements and generator trajectories are perfectly circular.}
	\label{fig.shu_rho_ie}
\end{figure}

\begin{figure}[!bp]  % 3 colonne
	\centering	
	\includegraphics[width=0.495\linewidth]{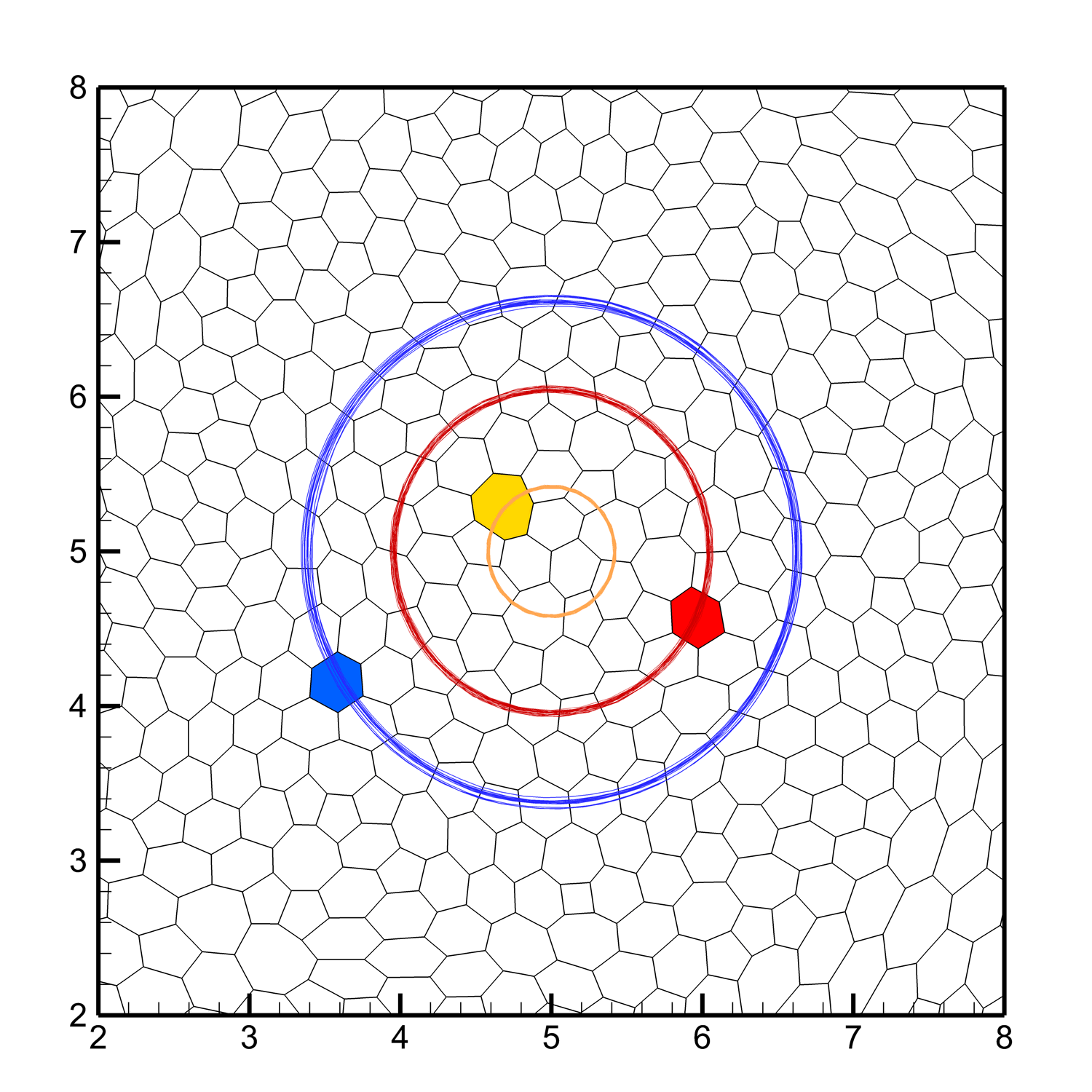}
	\includegraphics[width=0.495\linewidth]{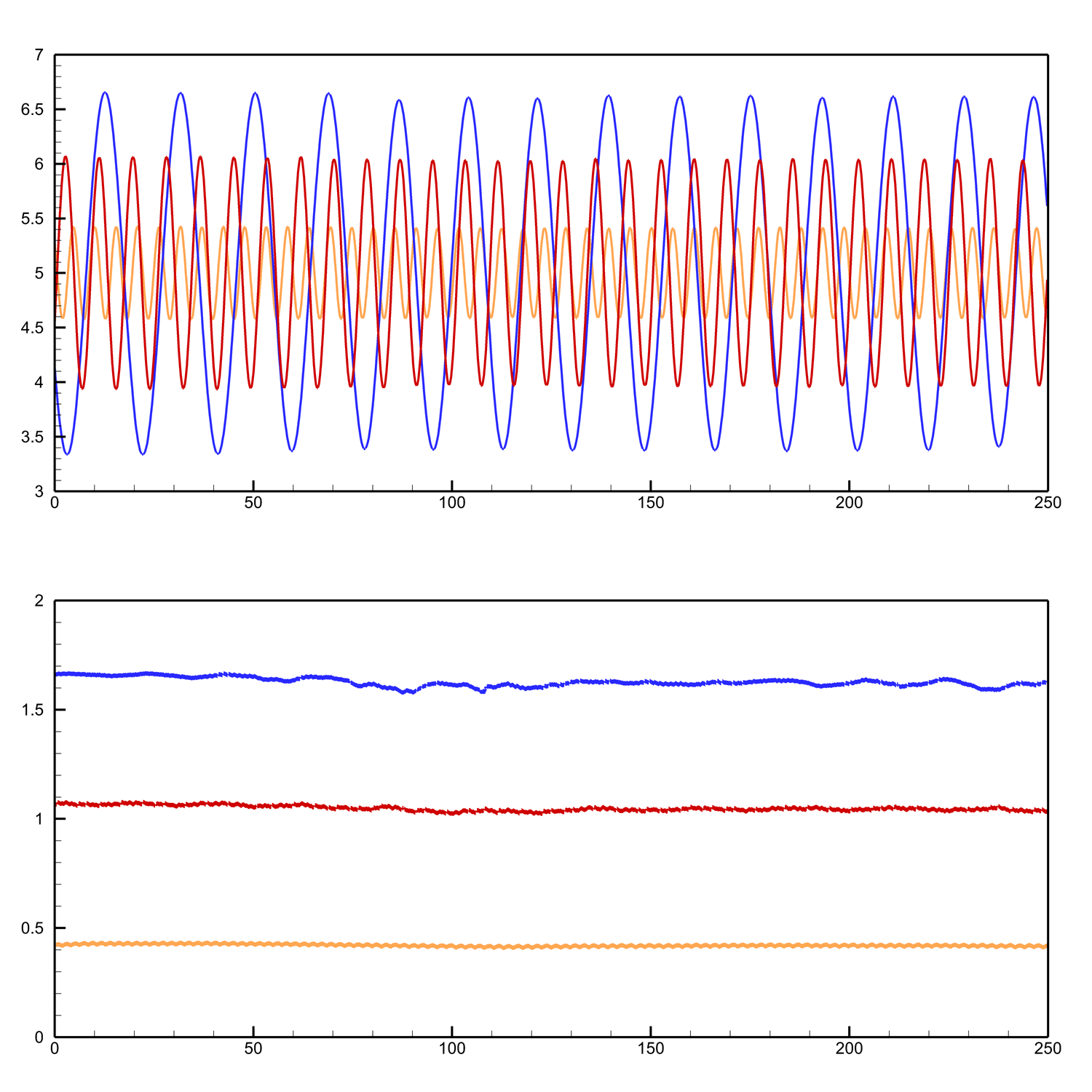}
	\caption{Stationary rotating vortex solved with our fourth order $P_3$ ALE-DG scheme on a moving Voronoi-type mesh of $957$ 
		elements with dynamical change of connectivity and with the  generators trajectories computed with fourth order of accuracy. 
		Left: We depict the trajectories (in Cartesian coordinates) of the generators of $3$ mesh elements 
		(those highlighted respectively in blue, violet and red) from time $t=0$ up to time $t=250$. During this 
		time interval the red mesh element completes $30$ revolutions about the origin. Right: we depict the $y$
		coordinates of the $3$ generators (top) and their radial coordinates (bottom). We would like to emphasize that the
		trajectories are circular (their radius is almost constant) for a very long evolution time.}
	\label{fig.trajectory}
\end{figure}

As a first test we consider a smooth isentropic vortex flow defined as similarly to~\cite{HuShuVortex1999}. 
The initial computational domain is the square $\Omega=[0;10]\times[0;10]$ and boundary conditions are of wall (slip) type everywhere.
The initial condition is given by some perturbations $\delta$ that are superimposed 
onto a homogeneous background state $\Q_0=(\rho,u,v,p)=(1,0,0,1)$, assuming that the entropy perturbation is zero, i.e. $\delta S= 0$. 
The perturbations for density and pressure are
\begin{equation}
	\label{rhopressDelta}
	\delta \rho = (1+\delta T)^{\frac{1}{\gamma-1}}-1, \quad \delta p = (1+\delta T)^{\frac{\gamma}{\gamma-1}}-1, 
\end{equation}
with the temperature fluctuation $\delta T = -\frac{(\gamma-1)\epsilon^2}{8\gamma\pi^2}e^{1-r^2}$ and the vortex strength $\epsilon=5$.
The velocity field is specified by
\begin{equation}
	\label{ShuVortDelta}
	\left(\begin{array}{c} \delta u \\ \delta v  \end{array}\right) = \frac{\epsilon}{2\pi}e^{\frac{1-r^2}{2}} \left(\begin{array}{c} -(y-5) \\ \phantom{-}(x-5)  \end{array}\right).
\end{equation}
This is a stationary equilibrium of the system so the exact solution coincides with the initial condition at any time.

Preserving this kind of vortical solution over long simulation times with 
minimal dissipation is a nontrivial task in a moving-mesh context. 
To achieve this result, we propose the use of a very high-order scheme (here an ADER-DG method of order~$4$) in a 
Lagrangian framework. We remark that the combination cannot be used with fixed topology, 
or advanced remapping techniques,
because the quality of the moving mesh subject to this constraint quickly deteriorates, 
as is clearly apparent in Figure~\ref{fig.shu_rho_ie_notopc}, 
where the simulation has to be stopped after about half a vortex rotation period.
This highlights the well-known fact that, for long time evolution, the mesh connectivity must be somehow updated.
In this work this is naturally achieved by means of space-time topology changes. 

Further, Figure~\ref{fig.shu_rho_ie} demonstrates that the treatment of topology changes via high-order 
integration over \textit{crazy} sliver elements is actually quite effective. Indeed one can note that 
the solution is visually the same at the beginning of the simulation and $500$ seconds after, even on a 
rather coarse mesh of only $957$ polygonal elements. 
Moreover, we take advantage of this test case to also emphasize the high precision of the mesh 
movement. The Voronoi-type polygonal cells, as well as the generator points, 
in fact can be observed to orbit along perfectly circular trajectories, as evidenced 
in Figures~\ref{fig.shu_rho_ie} and~\ref{fig.trajectory}.

\blue{
\subsubsection{Order of convergence}
\label{sec.order}

Finally, this stationary test case allows to show numerically the order of convergence of the proposed ALE-DG scheme with topology changes, reported in Table~\ref{tab.Euler_order} up to order $4$.
Furthermore, we present a quantitative comparison with the scheme applied in a purely Eulerian setting (i.e. on a fixed mesh) 
and with the classical direct ALE approach with fixed topology.
For the purpose of this test, we consider the domain $\Omega=[-10;30]\times[-10;30]$, 
covered with a Voronoi-type tessellation obtained as the centroid-based dual of a
Delaunay mesh generated by Ruppert's algorithm \cite{ruppert1993}. 
We report our results at time $t=4$ (the time at which the ALE simulations with fixed topology terminate due to mesh tangling) and $t=10$ 
(a long simulation time at which differences in mesh configuration become very significant).

It should be stressed that, 
due to the absence of discontinuous features or strong background flows, this test problem is not intended to highlight the capabilities
of moving mesh algorithms, but rather to show the high order convergence of 
the method on smooth flows, while highlighting the necessity for a changing mesh topology.

\begin{table}[!tp]
	\caption{ \blue{ Stationary vortex test case with final time $t=4$ and $t=10$. 
		We report here the order of convergence, on the variable $\rho$ in the $L_2$ norm, for our DG scheme 
		up to order 4 in the Eulerian case (left), for a standard Lagrangian method with fixed topology (middle) and for our ALE scheme with topology changes (right).
		In the last case we also report the total number of \textit{crazy} sliver elements that have been originated during the simulations:
		the high order of convergence is maintained also when many sliver elements appear in the mesh.
		At large times ($t=10$), the ALE scheme with topology changes, 
		produces numerical errors that are comparable or smaller than those obtained with the corresponding Eulerian method on a fixed mesh.
%Moreover, we highlight in bold which algorithm is giving the smallest absolute error at each time. We can notice that our Lagrangian scheme,
% at large times, i.e. after the mesh has gradually adapted to the flow field, shows smaller numerical errors w.r.t. the corresponding Eulerian method.
}
}
	\label{tab.Euler_order}   
    \centering
		\begin{tabular}{c|c|ccc|cc|ccccc}
			\hline\noalign{\smallskip}
			&&\multicolumn{3}{c}{Eulerian} & \multicolumn{2}{|c|}{ALE fixed}    &  \multicolumn{5}{c}{ALE+sliver} \\
			\noalign{\smallskip}\svhline\noalign{\smallskip}			
			&&\multicolumn{1}{c}{t=4} & \multicolumn{2}{c}{t=10}&\multicolumn{2}{|c|}{t=4} & \multicolumn{2}{c}{t=4} & \multicolumn{3}{c}{t=10} \\
		   	&$h$&$\epsilon(\rho)$ & $\epsilon(\rho)$ & $\mathcal{O}(L_2)$ & $\epsilon(\rho)$  & $\mathcal{O}(L_2)$ & \ \ Sliv & $\epsilon(\rho)$ &\ \ Sliv & $\epsilon(\rho)$ &$\mathcal{O}(L_2)$\\			
			\noalign{\smallskip}\hline\noalign{\smallskip}
			\multirow{2}{*}{ \rotatebox{90}{$P_1 \rightarrow \mathcal{O}2$} \, } 
			&\ \ 5.9E-1 \ \ & \ \ \BB{4.7E-2} \ \ &\ \  \BB{8.3E-2} &  -  \ \ &\ \  1.4E-1 &  -  \ \ &\ \ 111 & 4.8E-2      \ \ &\ \  293  & 9.0E-2      &  -  \\
			&\ \ 4.4E-1 \ \ & \ \ 2.8E-2      \ \ &\ \  4.6E-2      & 2.0 \ \ &\ \  1.1E-1 & 1.0 \ \ &\ \ 192 & \BB{2.4E-2} \ \ &\ \  514  & \BB{4.3E-2} & 2.6 \\		
			&\ \ 2.9E-1 \ \ & \ \ 1.0E-2      \ \ &\ \  1.8E-2      & 2.4 \ \ &\ \  4.6E-2 & 2.1 \ \ &\ \ 420 & \BB{9.3E-3} \ \ &\ \  1119 & \BB{1.5E-2} & 2.6 \\
			&\ \ 2.2E-1 \ \ & \ \ 4.9E-3      \ \ &\ \  8.0E-3      & 2.8 \ \ &\ \  2.4E-2 & 2.3 \ \ &\ \ 789 & \BB{4.6E-3} \ \ &\ \  2111 & \BB{6.5E-3} & 2.9 \\					
			\noalign{\smallskip}\hline\noalign{\smallskip}
			\multirow{2}{*}{ \rotatebox{90}{$P_2 \rightarrow \mathcal{O}3$} } 
			&\ \ 5.9E-1 \ \ & \ \ \BB{6.7E-3} \ \ &\ \  1.1E-2 &  -  \ \ &\ \  2.7E-2 &  -  \ \ &\ \ 97  & 7.3E-3 \ \ &\ \  277  & \BB{1.0E-2} &  -  \\
			&\ \ 4.4E-1 \ \ & \ \ \BB{2.8E-3} \ \ &\ \  3.9E-3 & 3.6 \ \ &\ \  1.3E-2 & 2.5 \ \ &\ \ 181 & 2.9E-3 \ \ &\ \  498  & \BB{3.2E-3} & 4.0 \\		
			&\ \ 2.9E-1 \ \ & \ \ \BB{9.6E-4} \ \ &\ \  1.1E-3 & 3.3 \ \ &\ \  6.8E-3 & 1.6 \ \ &\ \ 401 & 9.7E-4 \ \ &\ \  1066 & \BB{9.0E-4} & 3.2 \\
			&\ \ 2.2E-1 \ \ & \ \ \BB{4.0E-4} \ \ &\ \  4.1E-4 & 3.4 \ \ &\ \  3.3E-3 & 2.5 \ \ &\ \ 745 & 4.3E-4 \ \ &\ \  1981 & \BB{4.1E-4} & 2.8 \\
			\noalign{\smallskip}\hline\noalign{\smallskip}
			\multirow{2}{*}{ \rotatebox{90}{$P_3 \rightarrow \mathcal{O}4$} } 
			&\ \ 1.7E-0 \ \ & \ \ 5.6E-2      \ \ &\ \  8.0E-2 &  -  \ \ &\ \  5.8E-2 &  -  \ \ &\ \ 2  & \BB{5.3E-2} \ \ &\ \  4   & \BB{6.8E-2} &  -  \\
			&\ \ 1.1E-0 \ \ & \ \ \BB{1.2E-2} \ \ &\ \  2.2E-2 & 3.2 \ \ &\ \  2.8E-2 & 1.8 \ \ &\ \ 10 & 1.3E-2      \ \ &\ \  41  & \BB{1.9E-2} & 3.2 \\		
			&\ \ 8.7E-1 \ \ & \ \ \BB{4.7E-3} \ \ &\ \  6.6E-3 & 4.4 \ \ &\ \  2.0E-2 & 1.3 \ \ &\ \ 36 & 5.8E-3      \ \ &\ \  110 & \BB{6.1E-3} & 4.2 \\
			&\ \ 5.9E-1 \ \ & \ \ 1.1E-3      \ \ &\ \  1.3E-3 & 4.2 \ \ &\ \  5.4E-3 & 3.3 \ \ &\ \ 93 & \BB{1.1E-3} \ \ &\ \  257 & \BB{1.3E-3} & 3.9 \\
			\noalign{\smallskip}\hline\noalign{\smallskip}			
		\end{tabular}
\end{table}
\begin{figure}[!bp]
	\centering
	\includegraphics[width=0.33\linewidth]{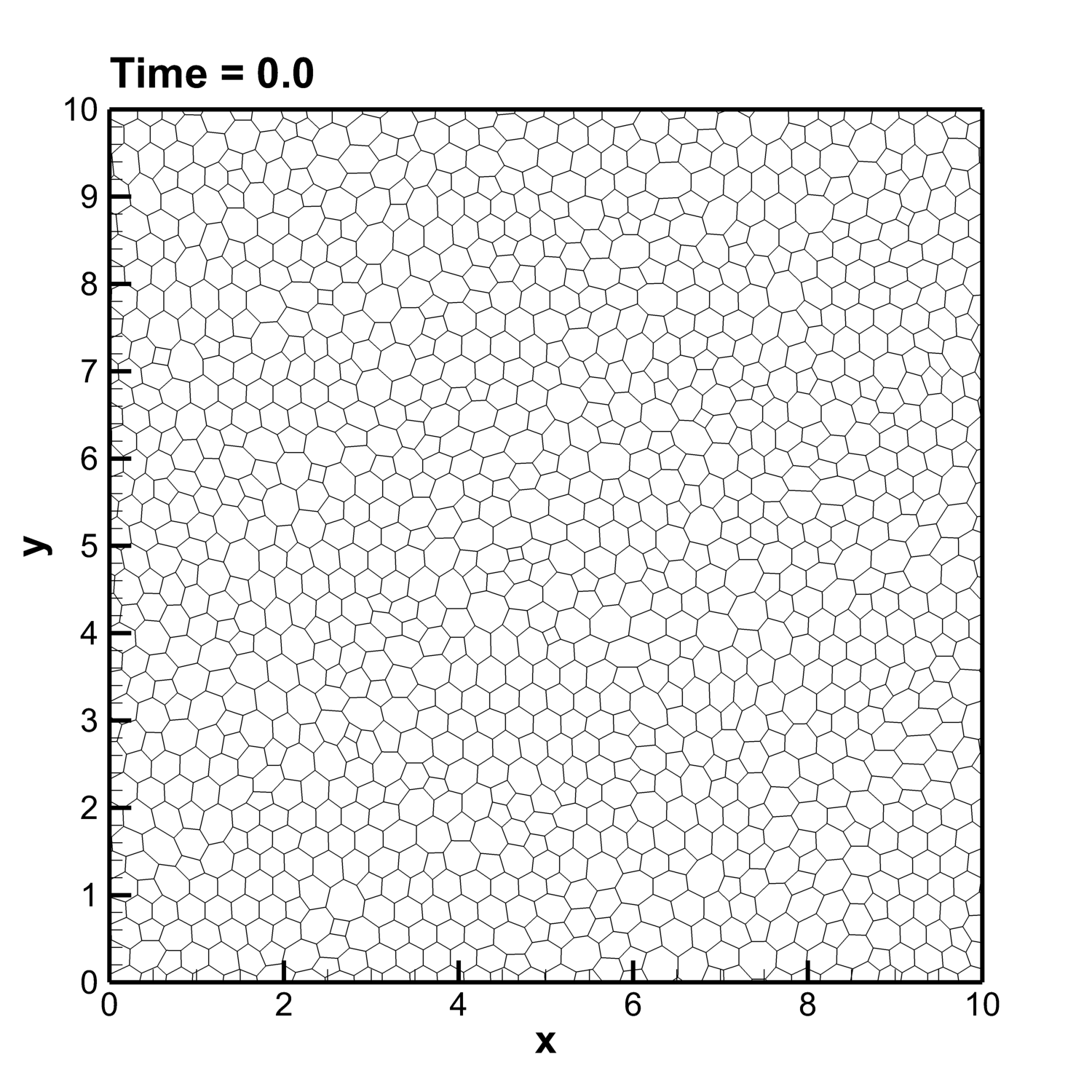}%
	\includegraphics[width=0.33\linewidth]{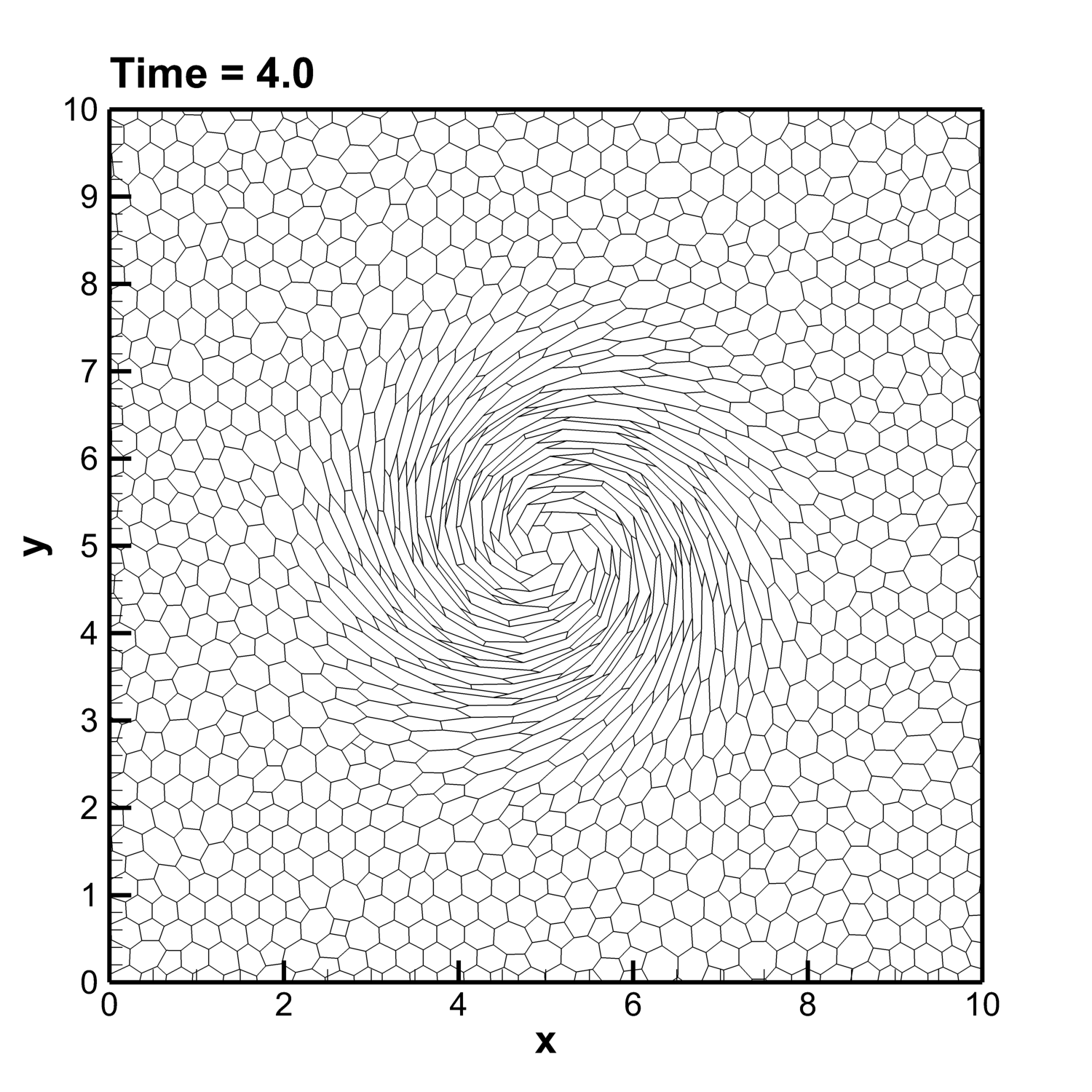}%
	\includegraphics[width=0.33\linewidth]{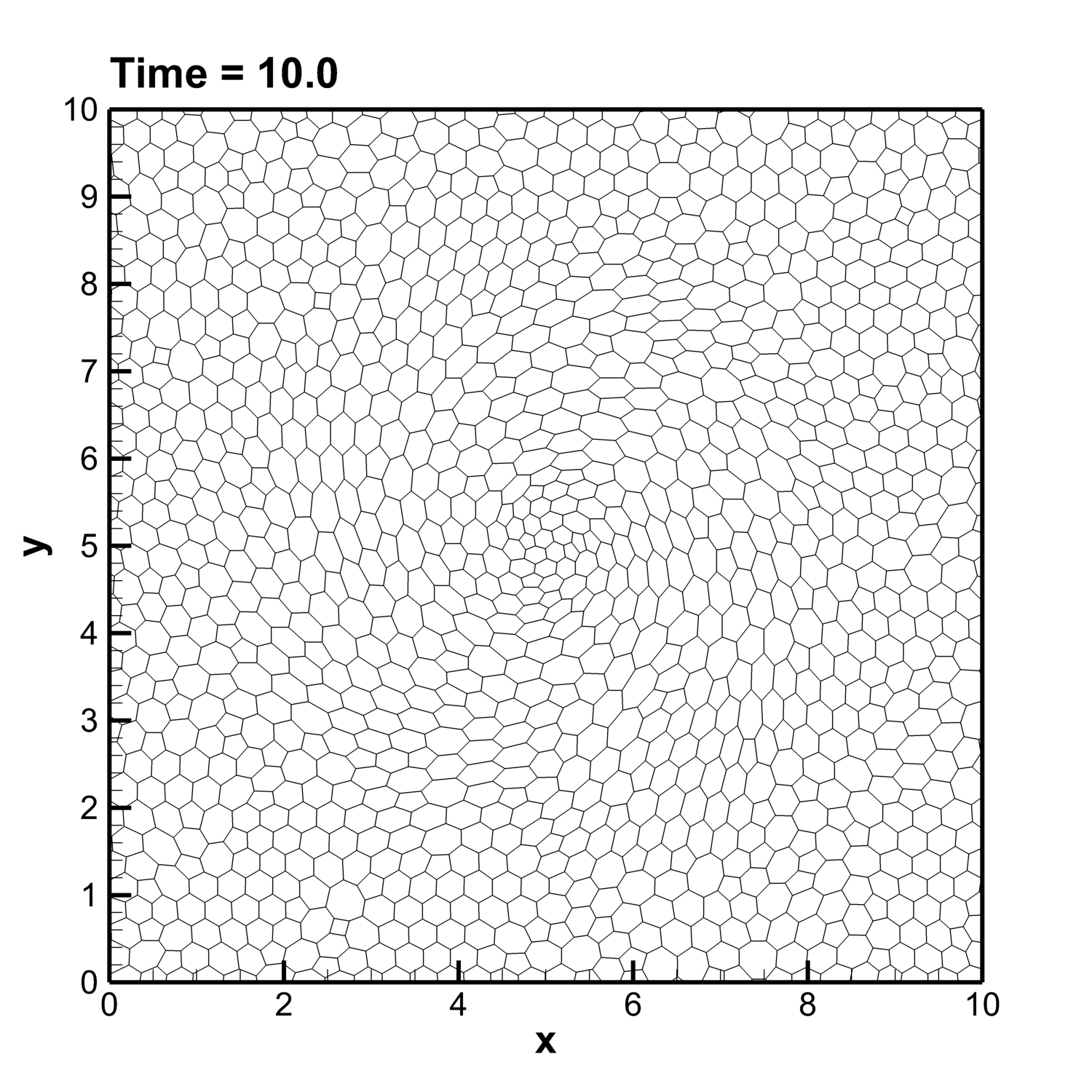}%
	\caption{\blue{Stationary vortex test case. We show here an example of a mesh employed for the 
	convergence test case of Section~\ref{sec.order} and of its evolution due to different Lagrangian schemes. 
	In particular, on the left we show the initial mesh, in the middle the mesh obtained with a standard direct 
	ALE scheme with fixed topology at time $t=4$, i.e. just before the simulation terminates due to mesh tangling, 
	and on the right the mesh obtained at time $t=10$ with our ALE algorithm dealing with topology changes, which has 
	gradually adapted to the fluid flow, optimizing the element shapes and allowing an increased precision for the DG scheme.}}
	\label{fig.convergencemesh}
\end{figure}

On the results, we note that the ALE method applied to a fixed mesh topology,
in addition to early termination around time $t=4$, as shown in Figure~\ref{fig.shu_rho_ie_notopc}, 
also suffers an increase in numerical errors, to the point that the correct order of convergence cannot be obtained when the mesh is
severely tangled.
Instead, the ALE algorithm presented in this work, 
not only deals with topology changes without accuracy losses, but 
in fact the mesh motion allows 
to gradually optimize the shape of the elements with respect to the flow field. This gradual optimization procedure,
translates into lower errors at large times with respect to the Eulerian scheme, 
for which the mesh is fixed to its initial generic configuration. 
We refer also to Figure~\ref{fig.convergencemesh} for a visual illustration of the different mesh motion approaches. 
%To emphasize which algorithm is giving the best result in each case, in Table~\ref{tab.Euler_order}, the smallest numerical errors (among the three schemes) at time $t=4$ and at time $t=10$ 
%have been highlighted in bold.

Finally, we would like to emphasize that in Table~\ref{tab.Euler_order} we show the numerical errors obtained at large computational times ($t=4$ and $t=10$), 
when computations have been carried out for thousands of timesteps and thousands of \textit{crazy} sliver elements have appeared (the total number is indicated in the Table), 
hence showing that the numerical method is genuinely high order accurate also when sliver elements are present.

}

\subsection{Sedov explosion problem}
\label{sec.sedov}

This test problem is a classic benchmark in the literature~\cite{LoubereSedov3D} and describes the evolution of a 
strong blast wave that is generated at the origin $\mathbf{O}=(x,y)=(0,0)$ of the computational 
domain $\Omega(0)=[0;1.2]\times[0;1.2]$. 
The difficulty of this benchmark is mainly due to the near zero pressure outer state
that may induce positivity-preservation problems. An exact solution based on self-similarity 
arguments is available from~\cite{Sedov}.

\begin{figure}[!bp]  % 3 colonne
	\centering
	\includegraphics[width=0.495\linewidth]{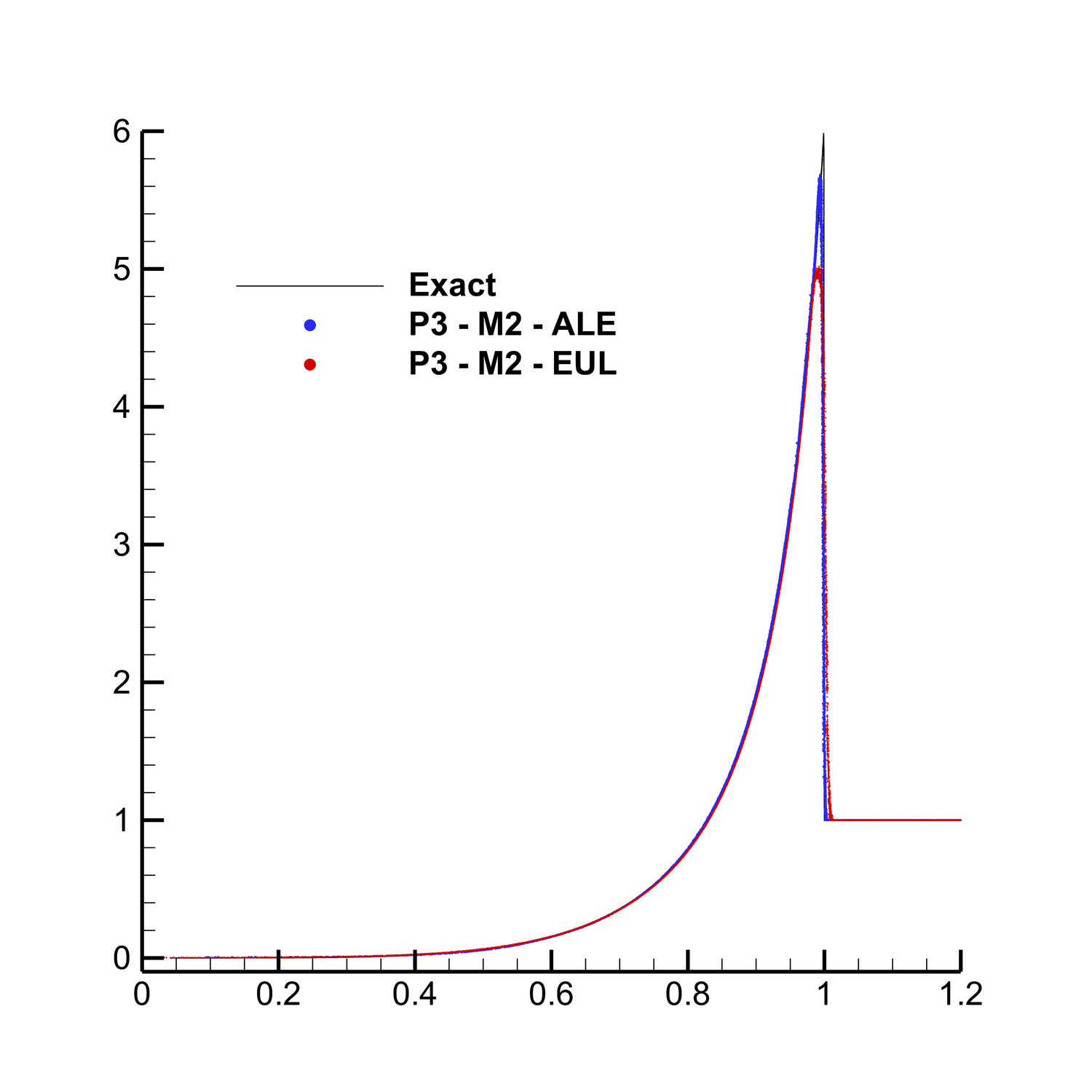}%		
	\includegraphics[width=0.495\linewidth]{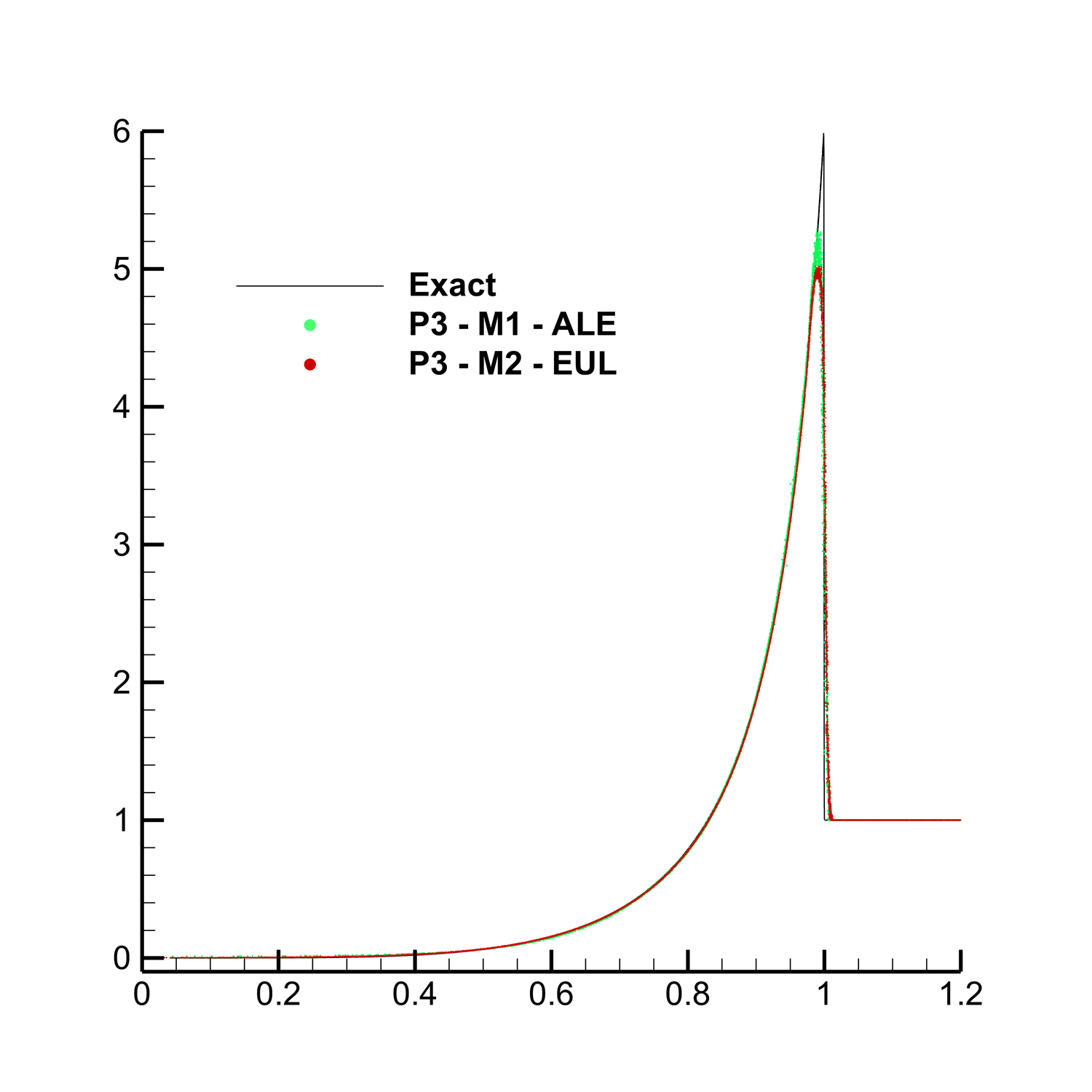}%		
	\caption{Sedov explosion problem. Comparison between the exact solution (black), the solution 
		obtained with a fourth order Eulerian $P_3$ DG method on the fine mesh $M2$  
		(red) and with our $P_3$ ALE-DG scheme both on $M2$ (blue) and $M1$ (green). 
		Our ALE scheme is more accurate than the Eulerian one even using coarser meshes.} 
	\label{fig.sedov_alevseul}
\end{figure}

The initial condition consists in a uniform density $\rho_0=1$ and a near zero pressure $p_0$ imposed 
everywhere except in the cell $V_{or}$ containing the origin $\mathbf{O}$ where it is given by
\begin{equation}
	p_{or} = (\gamma-1)\rho_0 \frac{E_{tot}}{|V_{or}|}, \quad \textnormal{ with } \
	E_{tot} =  0.979264,
	\label{eqn.p0.sedov}
\end{equation}
with $E_{tot}$ being the total energy concentrated in the cell containing the coordinate $\x=\mathbf{0}$.
We set $p_0=10^{-9}$ and solve this numerical test with a fourth order $P_3$ DG scheme; 
we employ a coarse mesh $M1$ made of $1345$ polygonal cells and a finer mesh $M2$ of $6017$ polygonal elements.

\begin{figure}[!bp]  % 3 colonne
	\centering
	\includegraphics[width=0.33\linewidth]{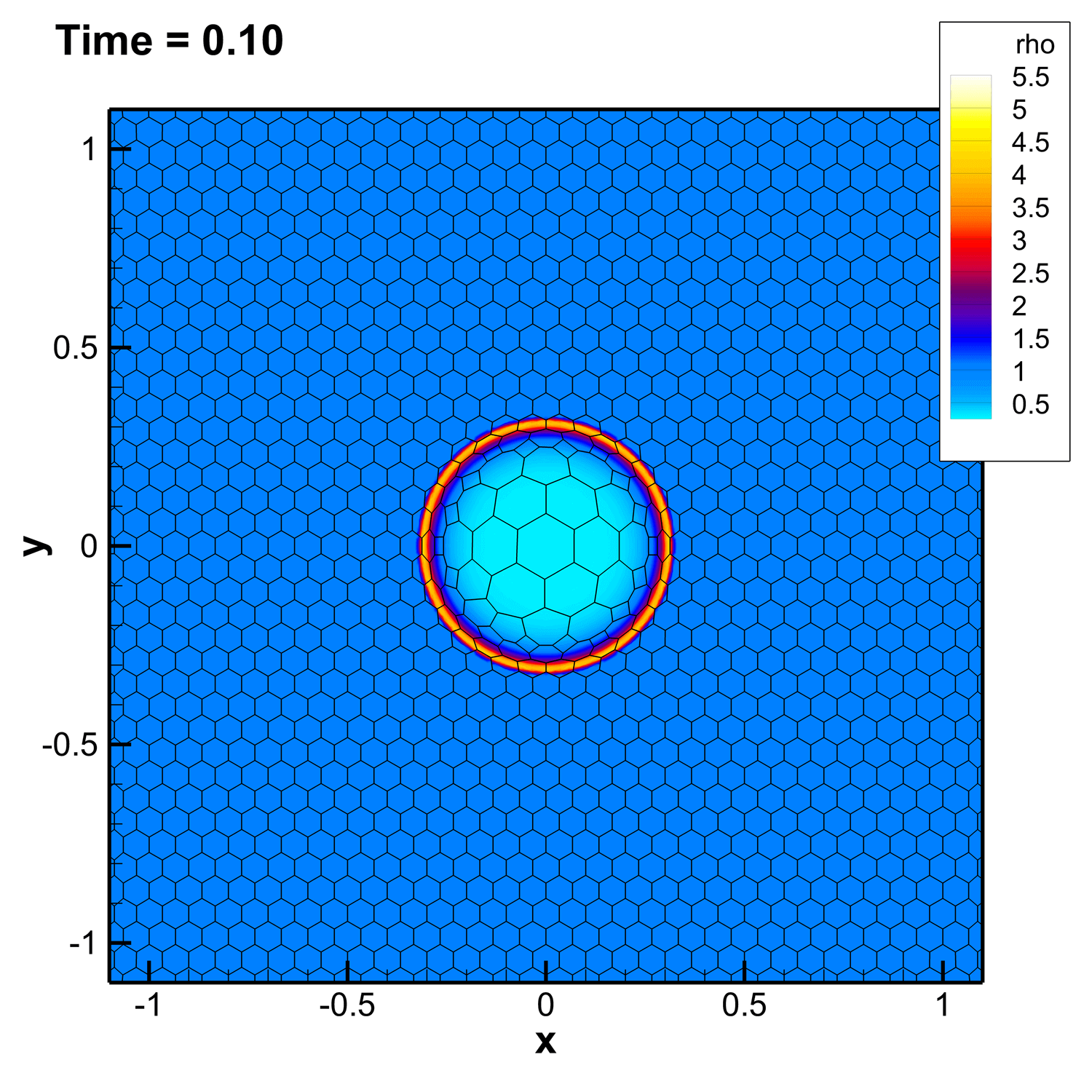}%
	\includegraphics[width=0.33\linewidth]{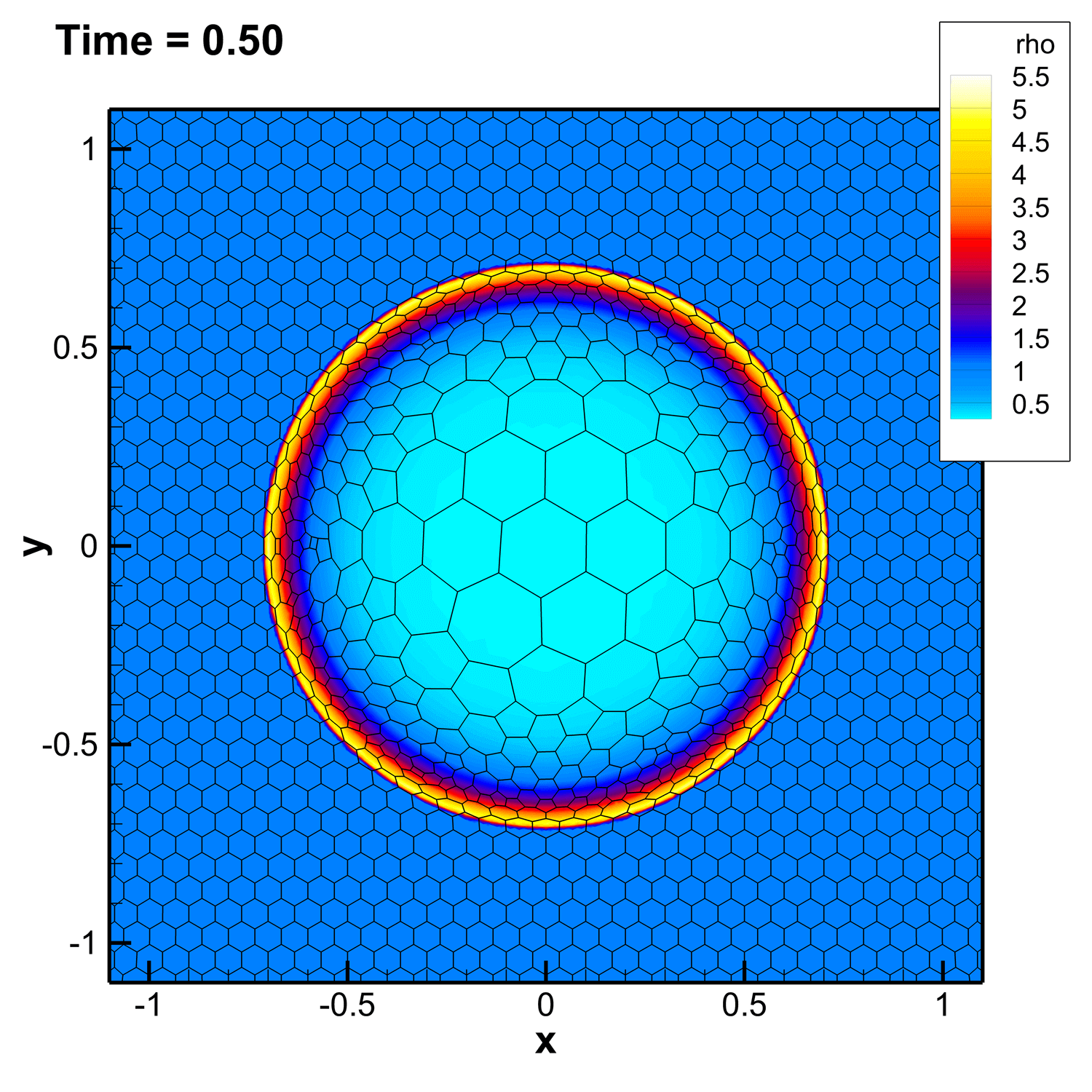}%
	\includegraphics[width=0.33\linewidth]{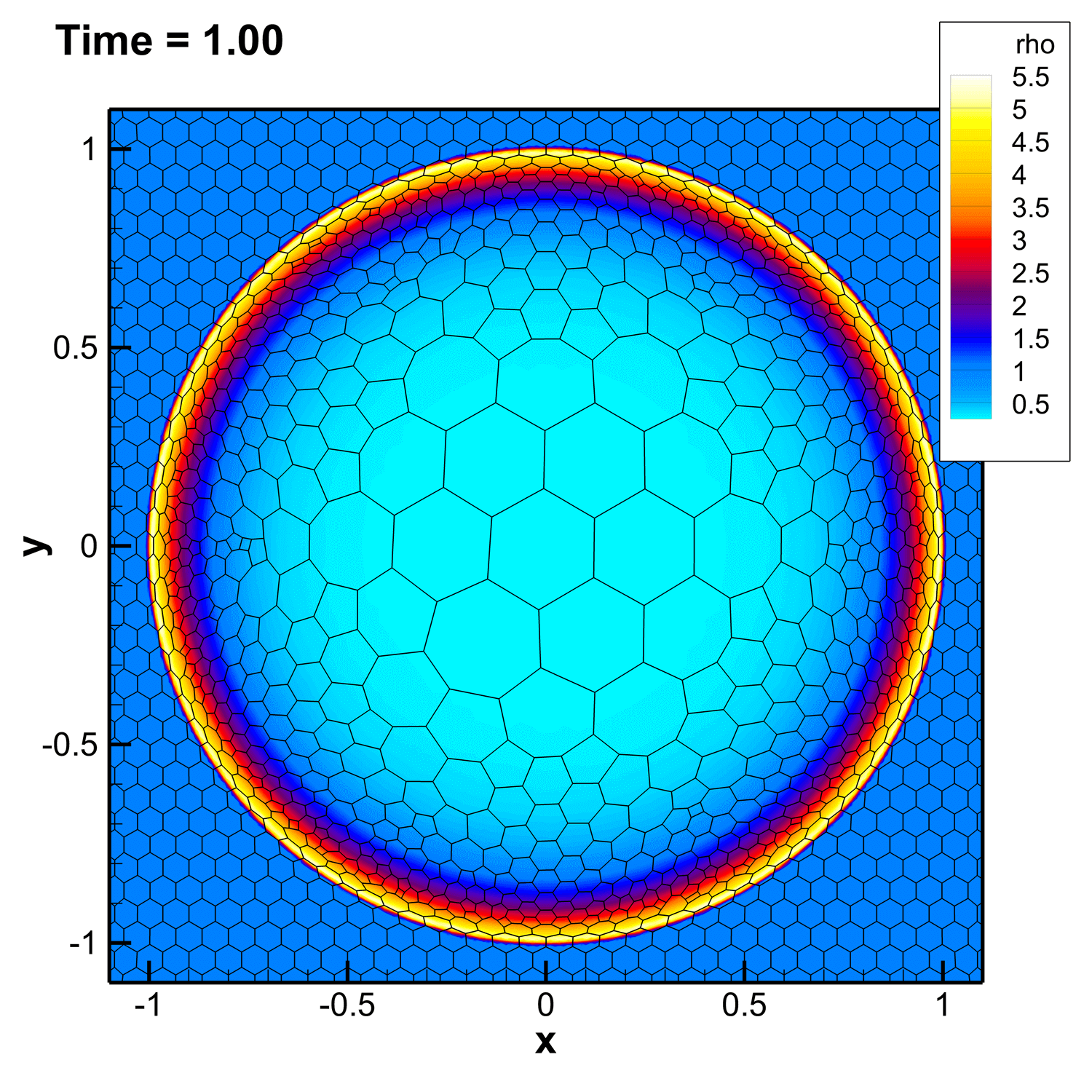}\\
	\includegraphics[width=0.33\linewidth]{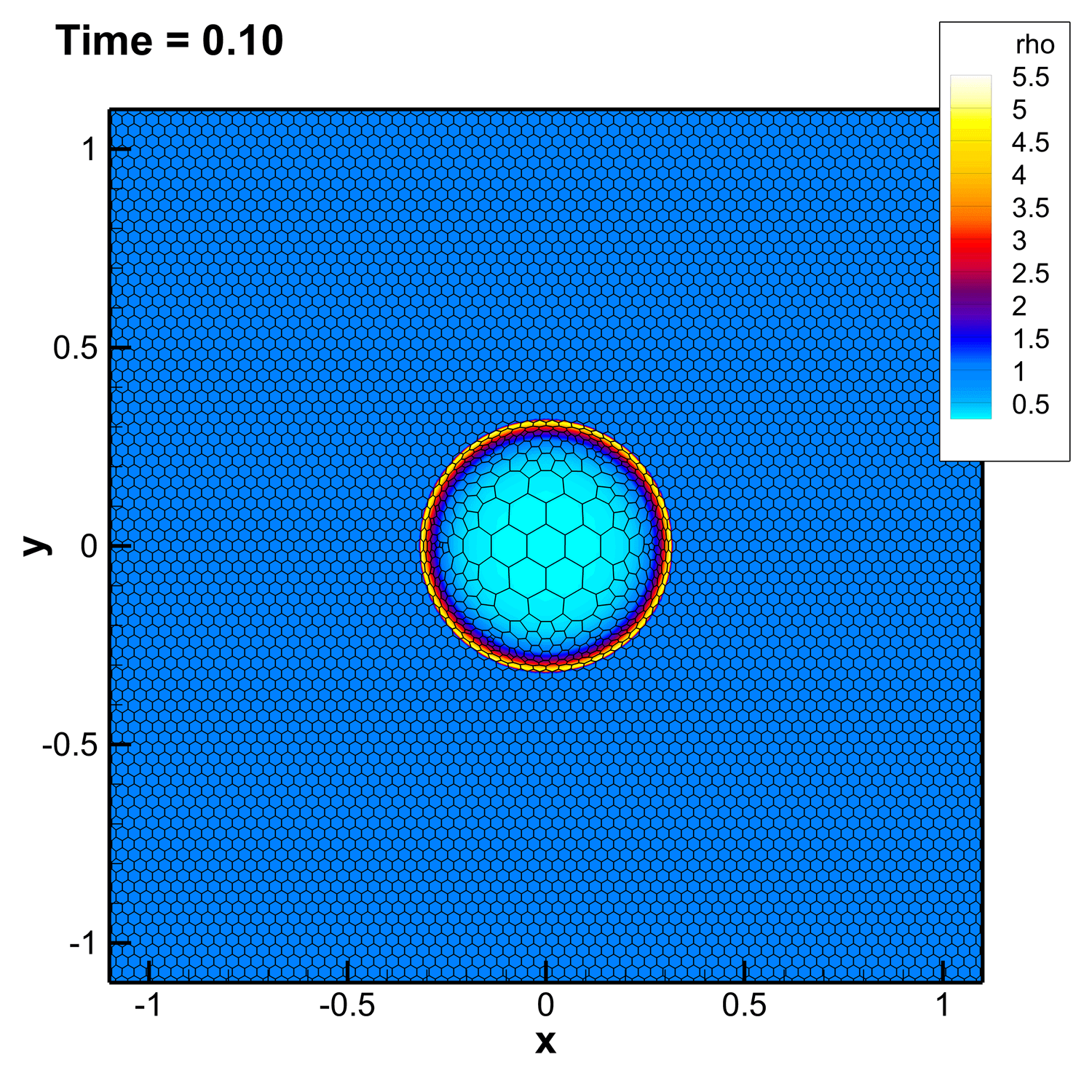}%
	\includegraphics[width=0.33\linewidth]{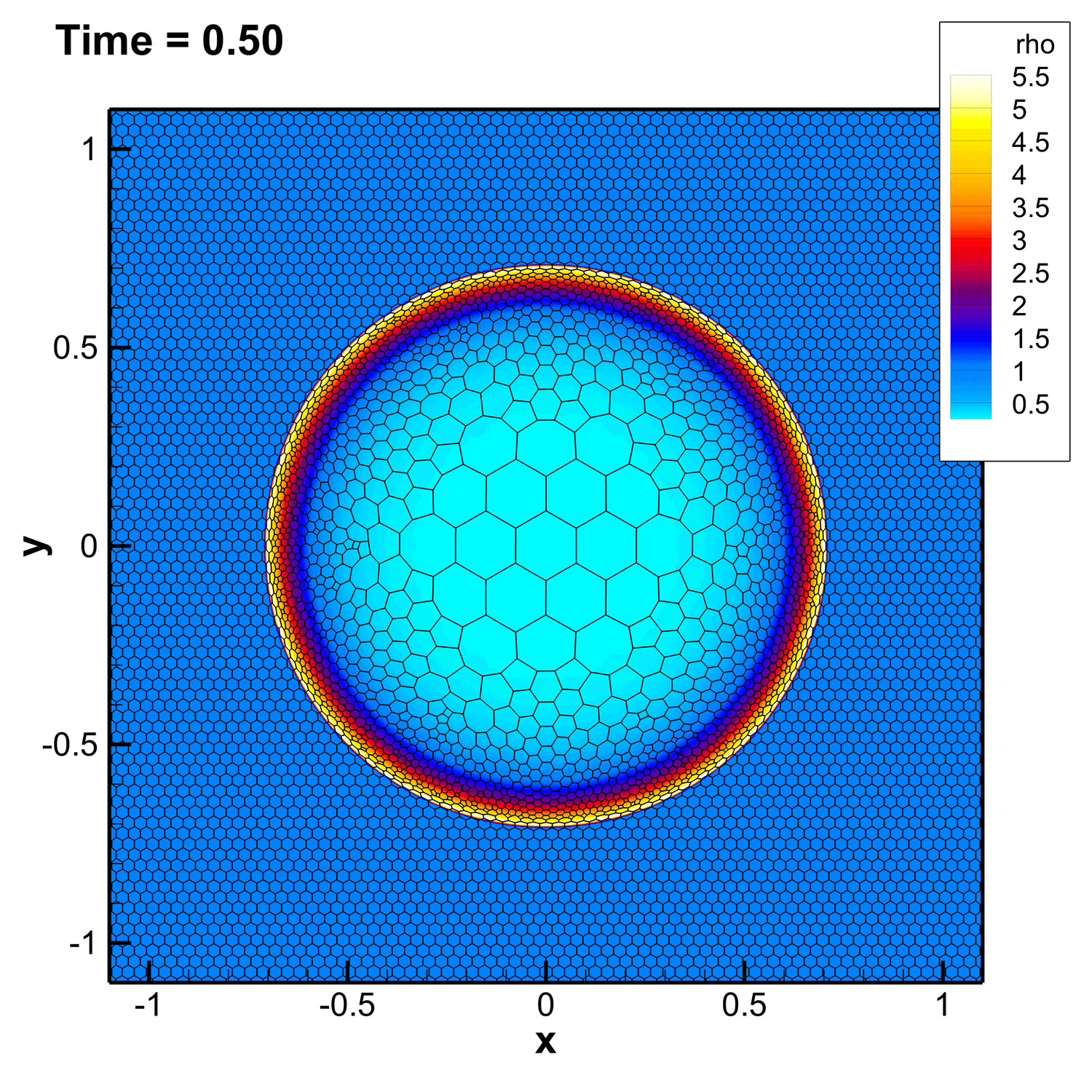}%
	\includegraphics[width=0.33\linewidth]{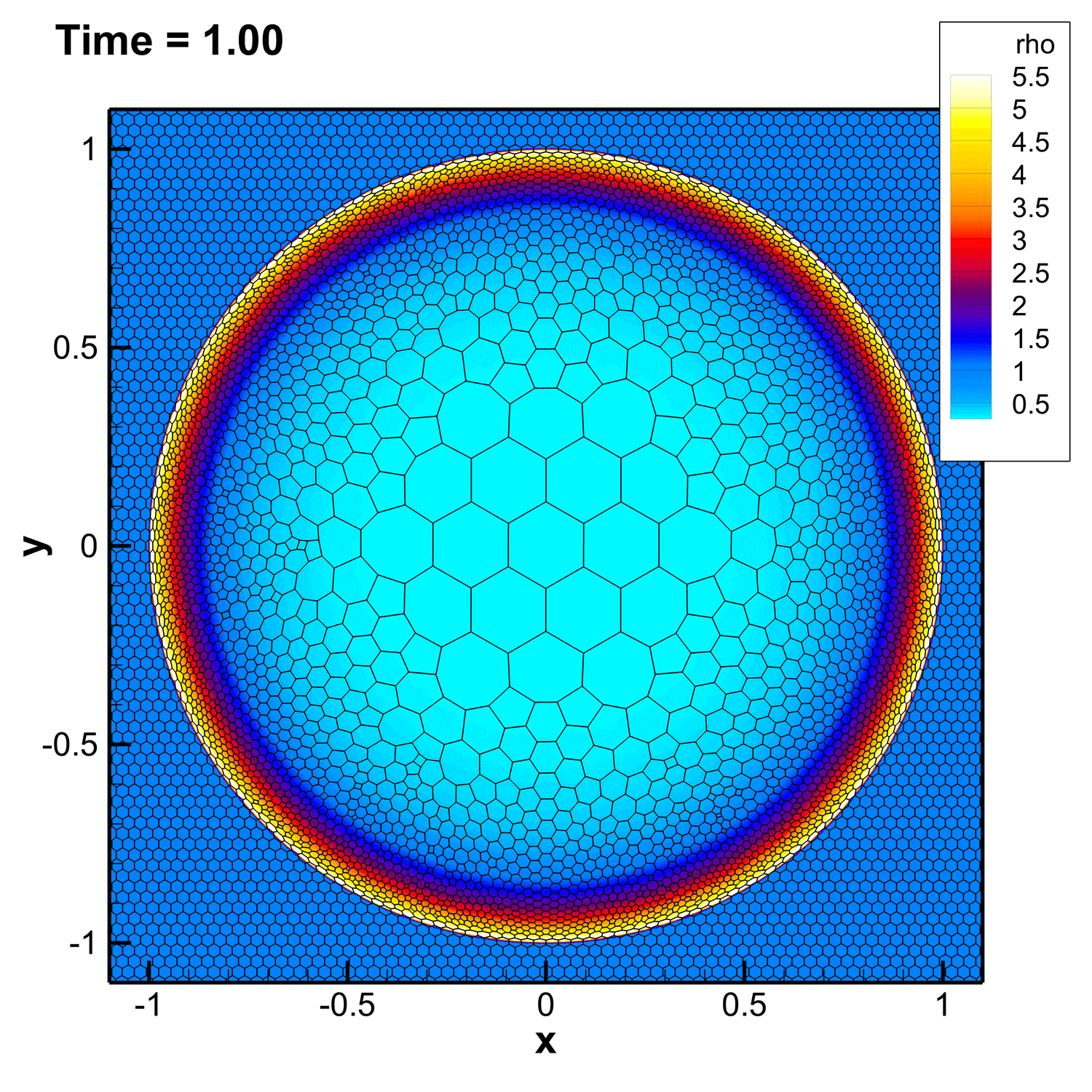}\\
	\caption{Sedov explosion problem. In this figure we show the density evolution and the corresponding 
		mesh movement at different output times computed with our $P_3$ ALE-DG scheme on the mesh $M1$ (top) and $M2$ (bottom).}
	\label{fig.sedov_evolution}
\end{figure}

\begin{figure}[!tp]  % 3 colonne
	\centering
	\includegraphics[width=0.495\linewidth]{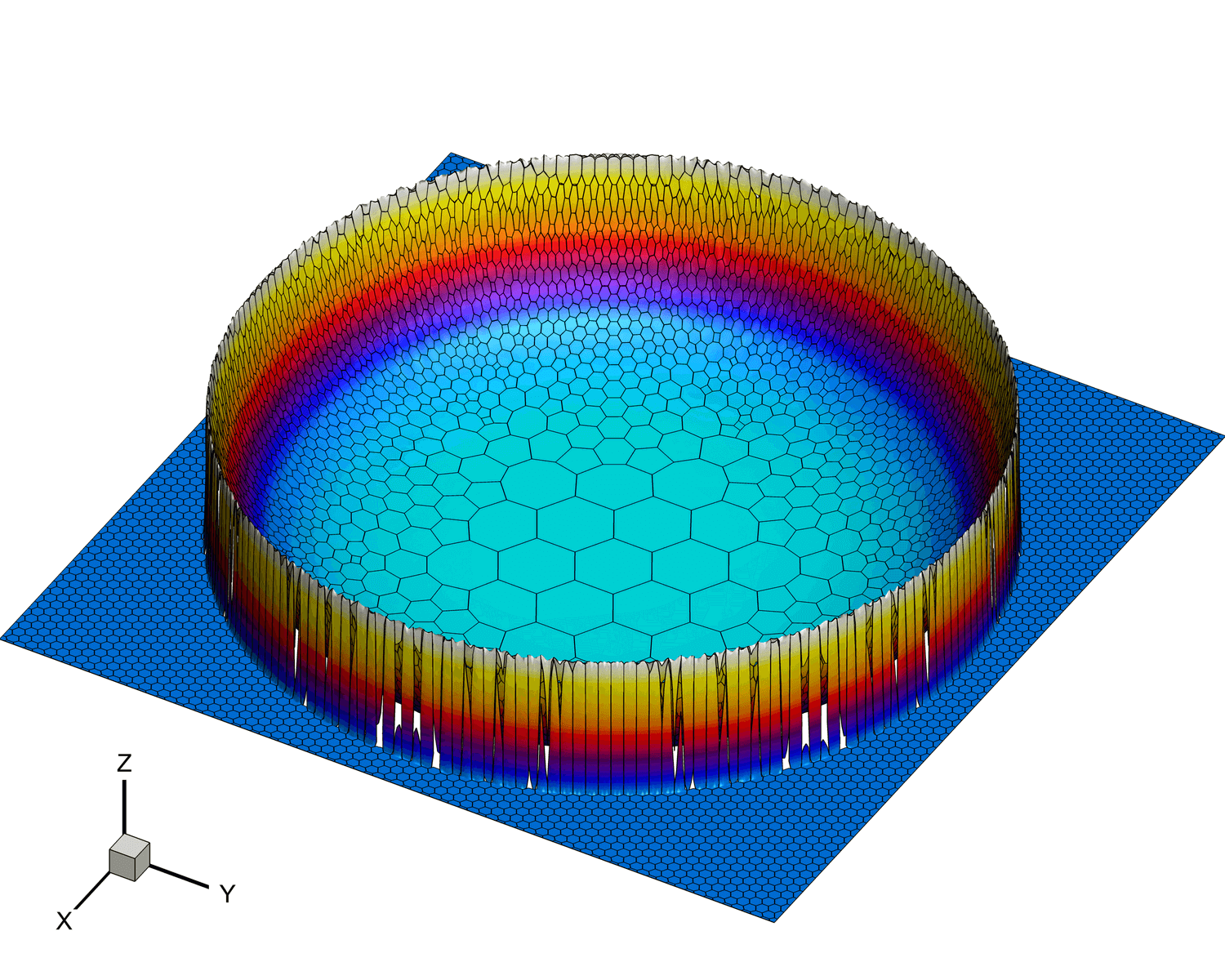}%		
	\includegraphics[width=0.495\linewidth]{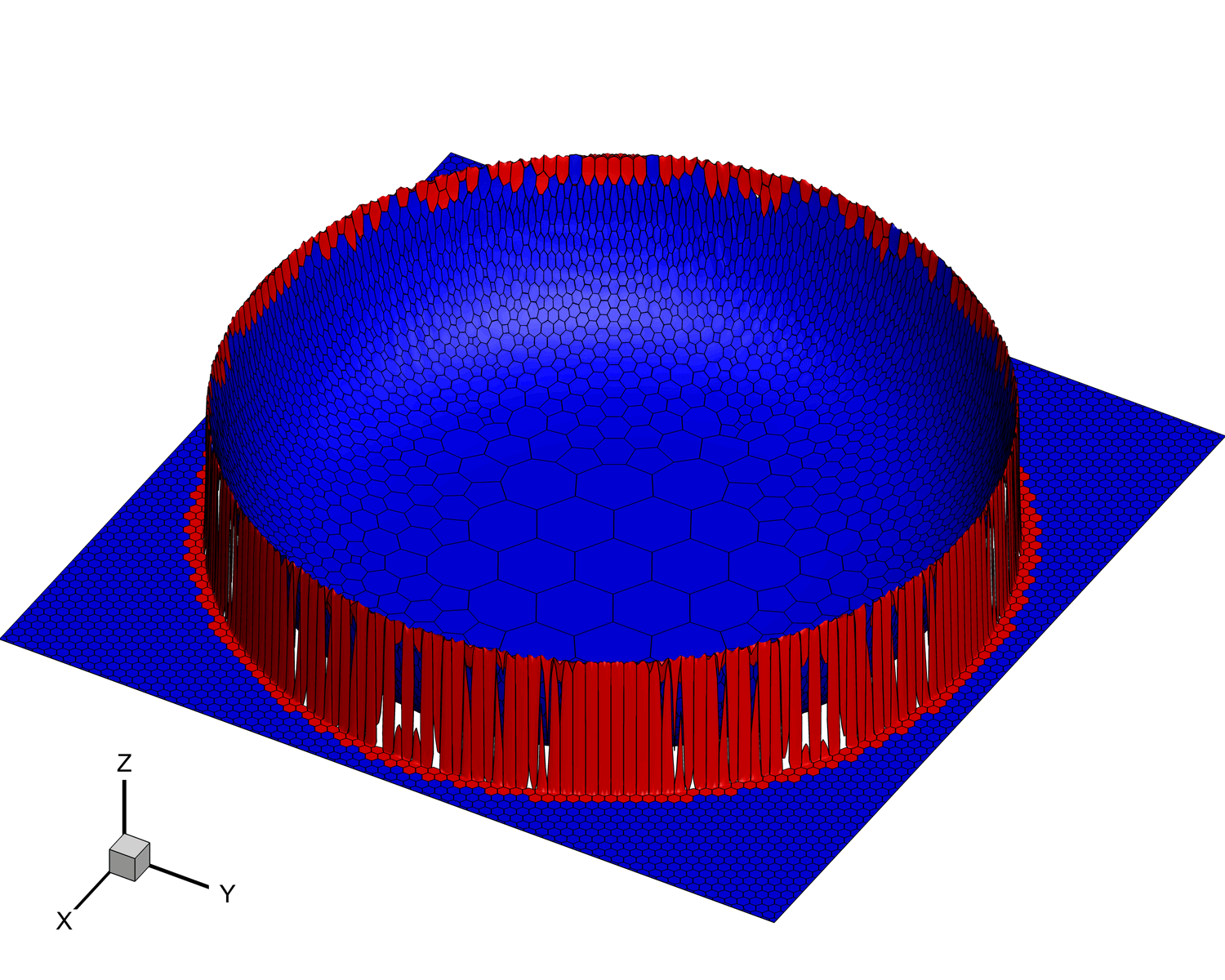}%	
	\caption{Sedov explosion problem. $3D$ density profile computed with our $P_3$ ALE-DG scheme 
		on $M2$. In particular, on the right, we highlight in red the so-called \textit{troubled} cells 
		marked by our detector on which the \textit{a posteriori} FV limiter has been employed.
		We make use of this image to further emphasize the robustness of our ALE schemes with topology 
		changes also in the presence of strong shock waves and near-zero pressure outer states.}
	\label{fig.sedov_3d}
\end{figure}

The density profiles are shown in Figure~\ref{fig.sedov_alevseul} for various output times $t~=~0.1, 0.5, 1.0$. 
The obtained results are in perfect agreement with the reference solution and the symmetry is very good despite using an unstructured grid, 
as opposed to a regular one built in polar coordinates. 
Also, one can note that the regularization procedure applied to the mesh elements does not compromise the natural 
expansion of the central cells expected in such an explosion problem.
Moreover, one can refer to Figure~\ref{fig.sedov_3d} for a comparison between our numerical solution (scatter plot) and the exact solution: 
the position of the shock wave and the density peak are perfectly captured. 

\begin{figure}[!bp]  % 3 colonne
	\centering
	\includegraphics[width=0.495\linewidth]{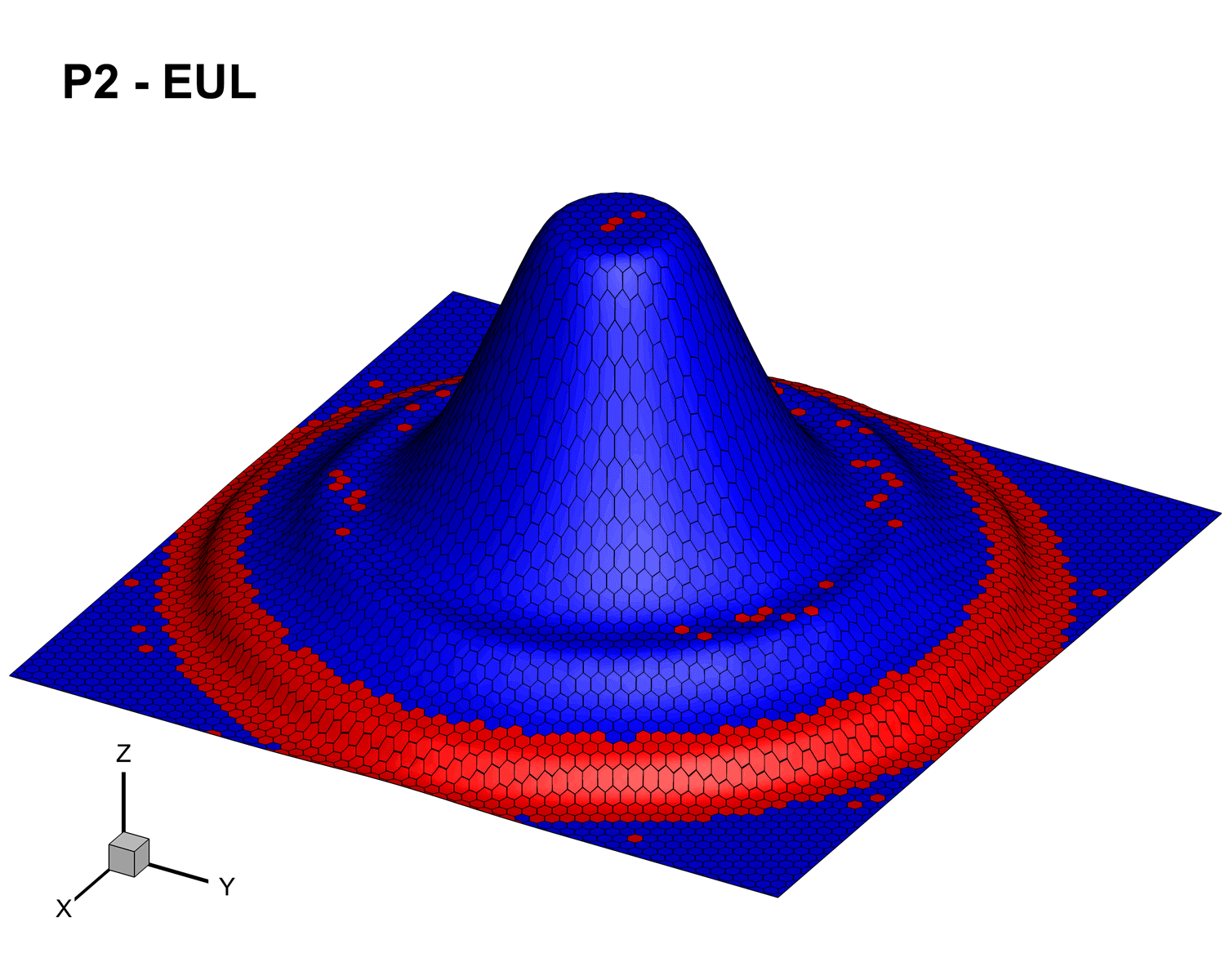}%
	\includegraphics[width=0.495\linewidth]{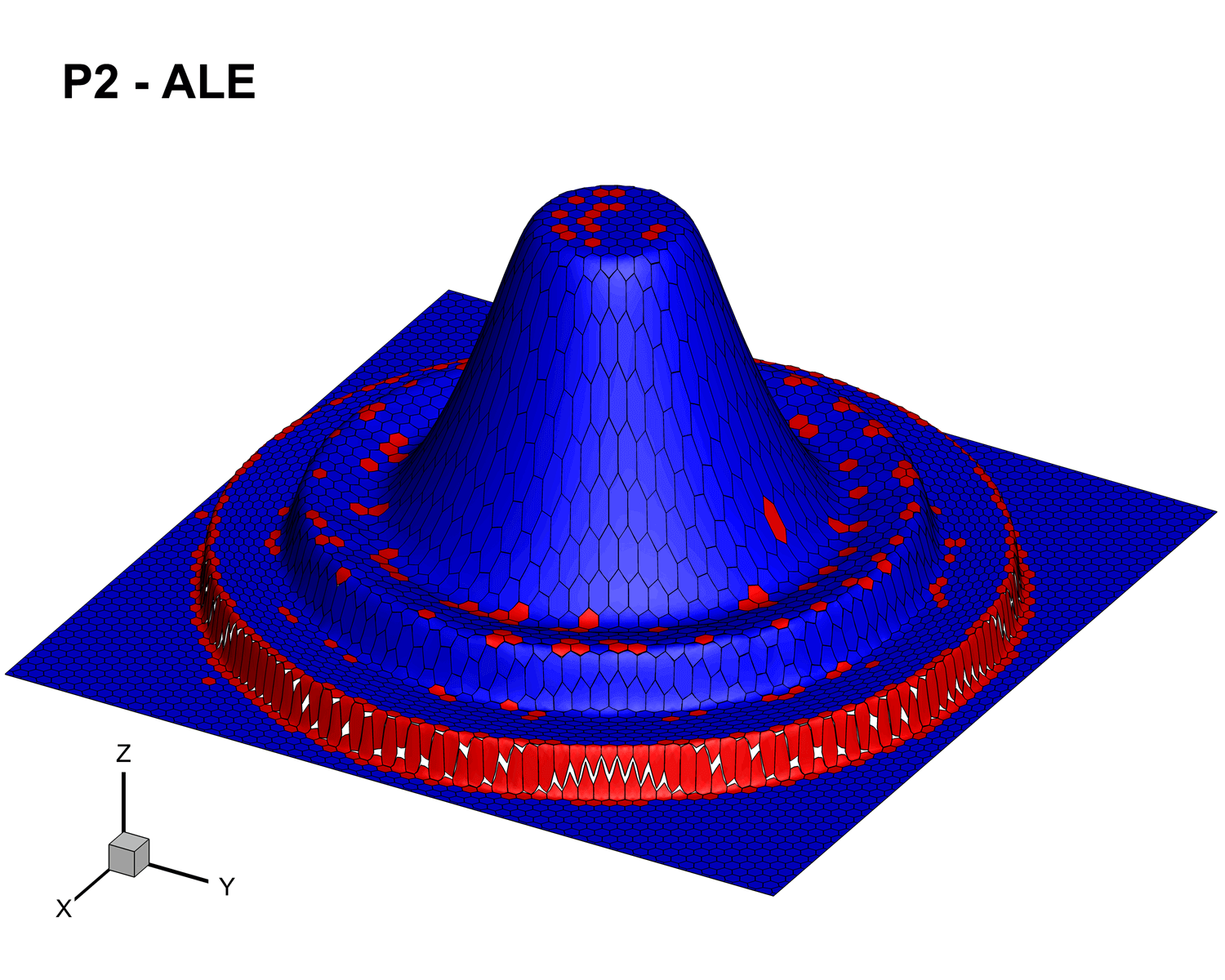}\\
	\includegraphics[width=0.495\linewidth]{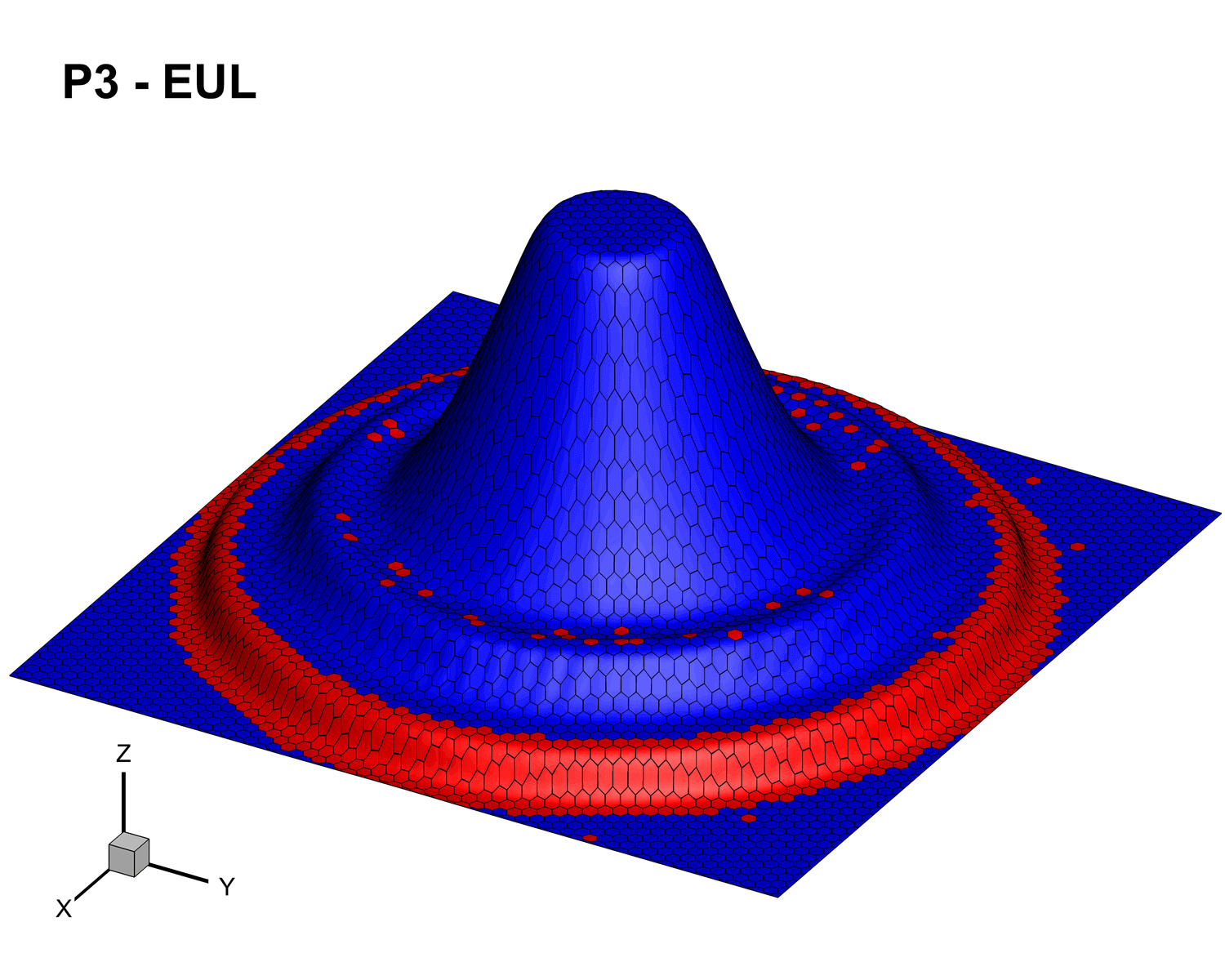}%
	\includegraphics[width=0.495\linewidth]{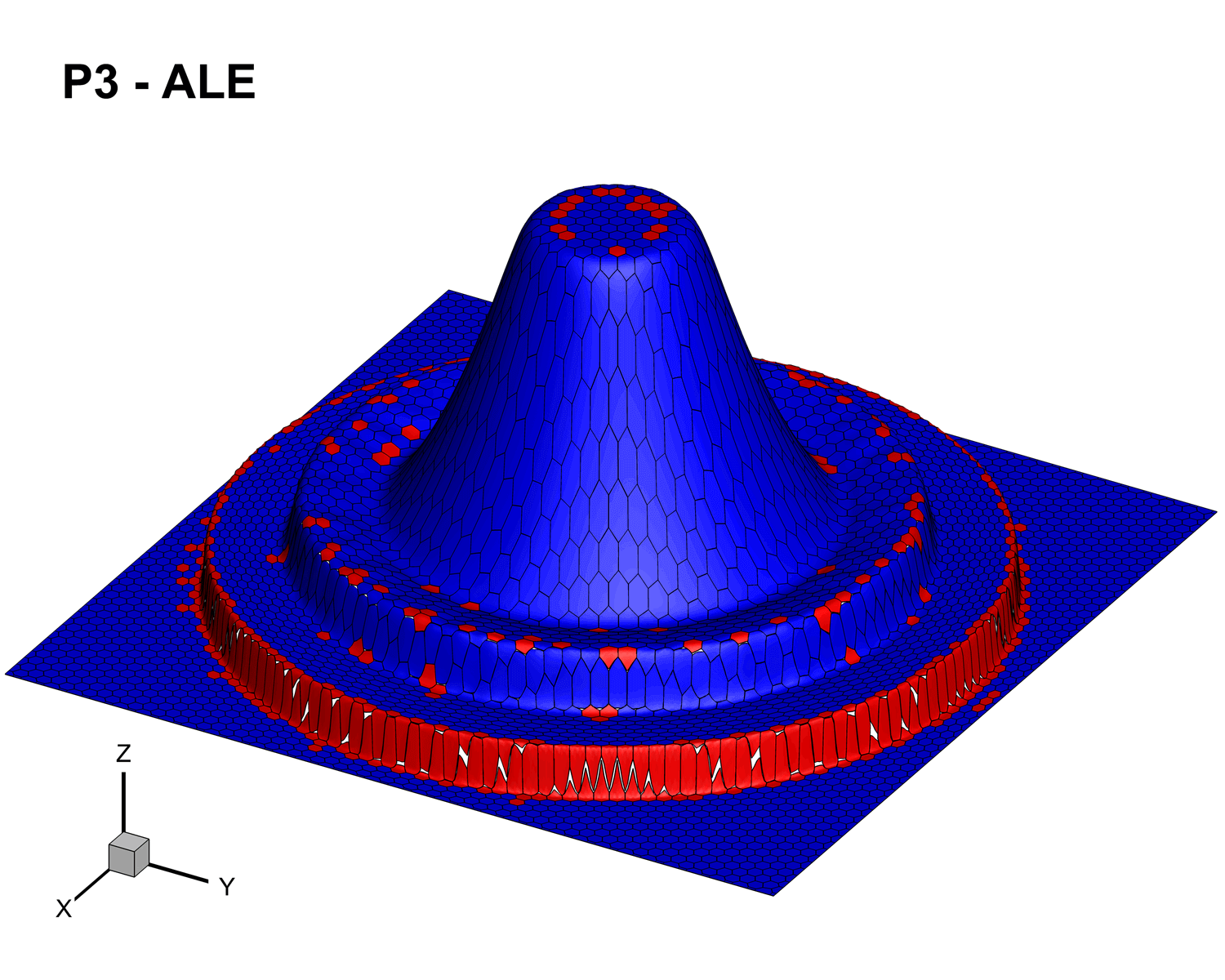}	
	\caption{Traveling Sod explosion problem. $3D$ density profile (z-axis) and limiter 
		activation (red cells), over a domain located in $[8.9;11.1]\times[-1.1;1.1]$ at the final time $t_f=0.25$, 
		obtained with $P_2$ and $P_3$ ADER-DG schemes run on a fixed Eulerian mesh (left) and our 
		direct ALE framework with topology changes (right). The difference on the numerical dissipation 
		between the Eulerian and the Lagrangian schemes is quite evident. We clarify that these 
		results are obtained solving the classical Sod explosion problem over a 
		high speed moving background.}
	\label{fig.sod_alevseul3D}
\end{figure}

In particular, we have chosen this test case in order to emphasize that Lagrangian schemes show a superior resolution
w.r.t. Eulerian ones even when both are compared at very high-order of accuracy, 
and furthermore that our direct ALE scheme results more accurate then the Eulerian method, 
even on a mesh ($M1$) that is coarser by a factor of two 
with respect to the finer mesh $M2$.

Finally, we refer to Figure~\ref{fig.sedov_3d} for the behavior of our \textit{a posteriori} sub-cell finite volume limiter, 
which activates only where the shock wave is located and is able to avoid any spurious oscillations or positivity problems, 
as can be noticed from the precise $3D$ density profile shown in Figure~\ref{fig.sedov_3d}.

\subsection{Traveling Sod-type explosion problem}
\label{sec.sod}

\begin{figure}[!bp]  
	\includegraphics[width=0.495\linewidth]{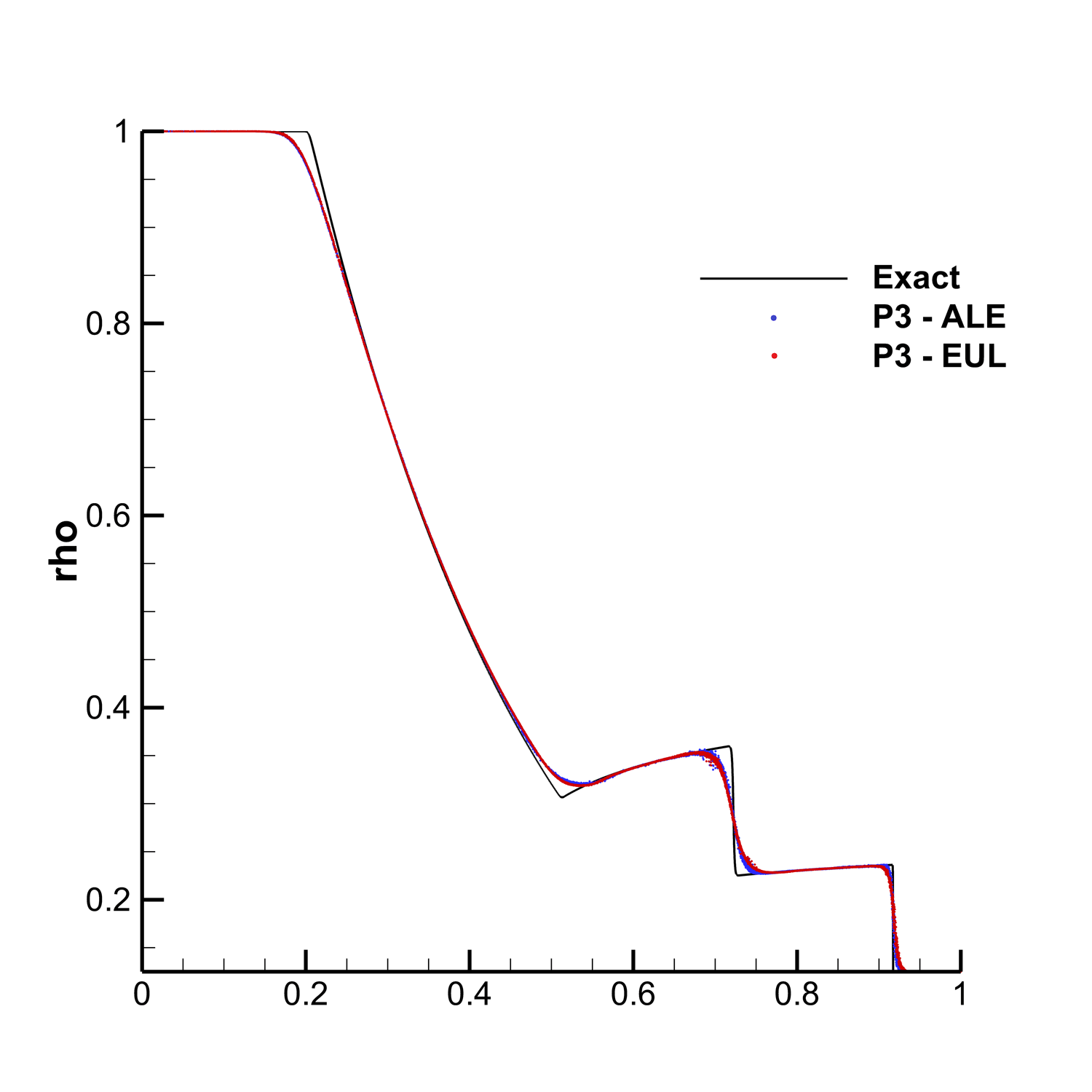}%	
	\includegraphics[width=0.495\linewidth]{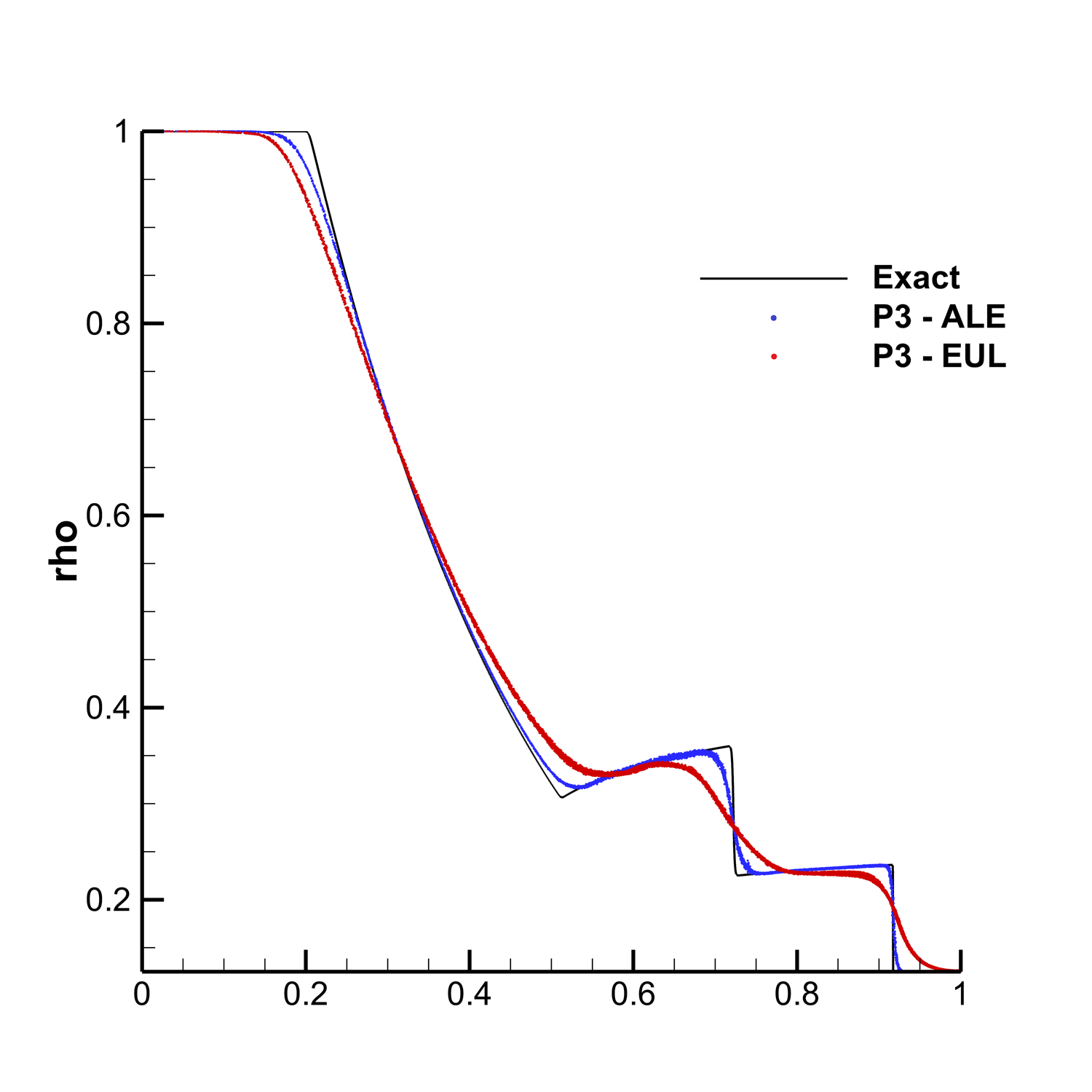}\\		
	\caption{Sod explosion problem: fixed background (left) and high speed traveling background (right).
		Comparison between the $P_3$ DG schemes on fixed Eulerian meshes (red) and in the moving-mesh ALE framework (blue). 
		This numerical results clearly explain that the Lagrangian schemes allow to obtain minimally
		dissipative results not only in a vanishing background flow, but even in a high speed one, and therefore
		that Lagrangian 
		methods discretely preserve the Galilean invariance of the equations. 
		On the contrary the influence of strong background flows on the solution obtained with Eulerian schemes is immediately apparent.}
	\label{fig.sod_alevseul}
\end{figure}

The explosion problem can be seen as a multidimensional extension of the classical Sod test case.
Here, we consider as computational domain a square of dimension $[-1.1;1.1]\times[-1.1;1.1]$ covered 
with a mesh made of $4105$ Voronoi-type elements, and the initial condition is composed of two 
different states, separated by a discontinuity at radius $r_d=0.5$ 
\begin{equation}
	\begin{cases} 
		\rho_L = 1, \quad \mathbf{u}_L, \quad p_L = 1, \quad  &\left\| \x \right\| \leq r_d \\
		\rho_R = 0.125, \quad \mathbf{u}_R, \quad  p_R = 0.1, \quad & \left\| \x \right\| > r_d .
	\end{cases}
\end{equation}
In addition, we aim at capturing the evolution of this explosion over a very high speed moving 
background (much higher than the speed of sounds): we impose $\mathbf{u}_L=\mathbf{u}_R=40$, so that 
at the final simulation time $t_f=0.25$ the square $[-1;1]\times[-1;1]$ will have been displaced by $5$ times its initial size.

We would like to underline that this test problem involves three different waves, therefore it allows 
each ingredient of our Lagrangian scheme to be properly checked. Indeed, we have \textit{(i)}~one cylindrical 
shock wave that is running towards the external boundary: our high-order scheme does not exhibit spurious
oscillations thanks to the \textit{a posteriori} sub-cell finite volume limiter;
\textit{(i)}~a rarefaction fan traveling in the opposite direction, which is well captured 
thanks to the high-order of accuracy of the DG scheme; and
\textit{(iii)}~an outward-moving contact wave, which is well resolved thanks to the Lagrangian nature 
of our scheme, in which, since the mesh moves together with the fluid flow, we can introduce a 
minimal dissipation when computing approximate Riemann fluxes. 

In addition, the high speed moving background allows to show the \textit{translational invariance} 
property of the Lagrangian schemes that indeed perfectly captures the three waves even when the 
explosion solution is moving at high speed, while the Eulerian scheme is heavily affected by the increased numerical dissipation.
Numerical evidence of the above statements can be found in Figure~\ref{fig.sod_alevseul3D}; 
moreover, in Figure~\ref{fig.sod_alevseul} we show that for this mild explosion 
it is really the background motion that requires the use of Lagrangian schemes, 
which, while still useful, would be instead not fundamental on a fixed background.

Finally, we want to remark that, despite the very high dissipation associated with the high base convective speed,
the overall symmetry of the solution, even in the Eulerian case, is not entirely compromised, 
thanks to the use of polygonal elements (see~\cite{boscheriAFE2022} for further discussion on the benefits of adopting polygonal meshes).

\section{Conclusion and outlook}
\label{sec.conclusion}

The accuracy of our results clearly show that the new combination of very high-order schemes with regenerated meshes, 
that allow topology changes, may open new perspectives in the fundamental research field of Lagrangian methods.

We would like to remark that the chosen simple test cases can be seen as prototypes of classical 
difficulties in astrophysical applications. Indeed, we have proposed here a method able to deal 
with long time simulation of vortical phenomena, as those necessary for the study of gas clouds evolving around black holes and neutron stars, 
and events, like explosions or interactions with near zero pressure states,
occurring in superposition with high speed background flows, as for supersonic or relativistic jets originating from 
proto-planetary nebulae, binary stars or nuclei of active galaxies. 

Future developments of this work will mainly concern the improvements of its robustness and effectiveness 
through mesh optimization and smoothing techniques~\cite{re2017interpolation,anderson2018high,morgan2018reducing,dobrev2020simulation,wang2015high}	and structure preserving algorithms~\cite{klingenberg2019arbitrary,Kapelli2014,castro2008well,dumbser2019glm,dumbser2020numerical,abgrall2022reinterpretation,gaburro2022high} so that future applications will effectively target in particular supersonic flows 
in aerodynamics~\cite{antoniadis2012high,tsoutsanis2015comparison} and 
astrophysics~\cite{gaburro2021well,Torsion2019,olivares2022new}, as well as 
fluid-structure interaction problems~\cite{basting2017extended,kikinzon2018establishing,durrwachter2021efficient}.

\begin{acknowledgement}
E.~Gaburro is member of the CARDAMOM team at the Inria center of the University of Bordeaux in France and S.~Chiocchetti is member of the INdAM GNCS group in Italy. 
E.~Gaburro gratefully acknowledges the support received from the European Union’s Horizon 2020 Research and Innovation Programme under the Marie Skłodowska-Curie Individual Fellowship \textit{SuPerMan}, grant agreement No. 101025563. S.~Chiocchetti acknowledges the support obtained by the Deutsche Forschungsgemeinschaft (DFG) via the project DROPIT, grant no. GRK 2160/2.
\end{acknowledgement}

\bibliographystyle{plain}
\bibliography{references}

\end{document}